# SIEVE EMPIRICAL LIKELIHOOD RATIO TESTS FOR NONPARAMETRIC FUNCTIONS

By Jianqing Fan[1] and Jian Zhang[2]

*Princeton University and Chinese University of Hong Kong, and University of Kent at Canterbury, EURANDOM and Chinese Academy of Sciences*

Generalized likelihood ratio statistics have been proposed in Fan, Zhang and Zhang [*Ann. Statist.* **29** (2001) 153–193] as a generally applicable method for testing nonparametric hypotheses about nonparametric functions. The likelihood ratio statistics are constructed based on the assumption that the distributions of stochastic errors are in a certain parametric family. We extend their work to the case where the error distribution is completely unspecified via newly proposed sieve empirical likelihood ratio (SELR) tests. The approach is also applied to test conditional estimating equations on the distributions of stochastic errors. It is shown that the proposed SELR statistics follow asymptotically rescaled $\chi^2$-distributions, with the scale constants and the degrees of freedom being independent of the nuisance parameters. This demonstrates that the Wilks phenomenon observed in Fan, Zhang and Zhang [*Ann. Statist.* **29** (2001) 153–193] continues to hold under more relaxed models and a larger class of techniques. The asymptotic power of the proposed test is also derived, which achieves the optimal rate for nonparametric hypothesis testing. The proposed approach has two advantages over the generalized likelihood ratio method: it requires one only to specify some conditional estimating equations rather than the entire distribution of the stochastic error, and the procedure adapts automatically to the unknown error distribution including heteroscedasticity. A simulation study is conducted to evaluate our proposed procedure empirically.

Received September 2001; revised September 2003.

[1]Supported in part by NSF Grants DMS-01-96041 and DMS-02-04329 and the RGC Grant CUHK-4299/00P of HKSAR.

[2]Supported in part by the National Natural Science Foundation of China and a grant from the research programme at EURANDOM.

**1. Introduction.** Over the last two decades, many computationally intensive nonparametric techniques and theories have been boldly developed to exploit possible hidden structures and to reduce modeling biases of traditional parametric methods. Methods such as local polynomial fitting, spline approximations and orthogonal series expansions as well as dimensionality reduction techniques have been studied in great depth in various statistical contexts. Yet there are no generally applicable methods available for the inferences in nonparametric models. Various efforts have been made in the literature on nonparametric hypothesis testing. See, for example, Bickel and Ritov (1992), Eubank and Hart (1992), Härdle and Mammen (1993), Azzalini and Bowman (1993), Fan (1996), Fan and Li (1996), Spokoiny (1996), Inglot and Ledwina (1996), Kallenberg and Ledwina (1997) and Horowitz and Spokoiny (2001, 2002), among others. For an overview, see the recent book by Hart (1997). Adaptive minimax rate results are obtained by various authors, including Fan (1996), Spokoiny (1996), Horowitz and Spokoiny (2001), Fan and Huang (2001) and Fan, Zhang and Zhang (2001). However, most of the studies focus only on the one-dimensional nonparametric regression problem. They are difficult to extend to multivariate semiparametric and nonparametric models.

In an effort to derive a generally applicable testing procedure for multivariate semiparametric and nonparametric models, Fan, Zhang and Zhang (2001) proposed generalized likelihood ratio tests. The work is motivated by the fact that the nonparametric maximum likelihood ratio test may not exist in many nonparametric problems. Further, even if it exists, it is not optimal even in the simplest nonparametric regression setting. Generalized likelihood ratio statistics, obtained by replacing unknown functions by reasonable nonparametric estimators, rather than the MLE as in the parametric setting, have several nice properties. In the varying coefficient model

$$Y = a_1(U)X_1 + \cdots + a_p(U)X_p + \varepsilon, \tag{1.1}$$

where $(U, X_1, \ldots, X_p)$ are independent variables and $Y$ is the response variable, Fan, Zhang and Zhang (2001) unveil the following Wilks phenomenon: The asymptotic null distributions are independent of nuisance functions and follow a $\chi^2$-distribution (in a generalized sense) for testing the homogeneity

$$H_0 : a_1(\cdot) = \theta_1, \ldots, a_p(\cdot) = \theta_p \tag{1.2}$$

and for testing the significance of variables, such as

$$H_0 : a_1(\cdot) = a_2(\cdot) = 0. \tag{1.3}$$

In other words, the generalized likelihood ratio statistic $\lambda_n$ follows asymptotically a rescaled $\chi^2$-distribution in the sense that $(2b_n)^{-1/2}(r_K\lambda_n - b_n) \xrightarrow{\mathcal{L}} N(0,1)$ for a sequence $b_n \to \infty$ and a constant $r_K$. We will use the notation $r_K \lambda_n \overset{a}{\sim} \chi^2_{b_n}$ to denote the result. The significance of the result is that the scale constant $r_K$ and the degrees of freedom $b_n$ are independent of nuisance parameters, such as the joint density of $(U, X_1, \ldots, X_p)$ and the parameters $\theta_1, \ldots, \theta_p$ in (1.2) and the functions $a_3(\cdot), \ldots, a_p(\cdot)$ in (1.3). This Wilks phenomenon is the key to the success of the classical maximum likelihood ratio tests for parametric problems. With the above newly discovered Wilks phenomenon in nonparametric models, the $P$-values can easily be computed by using either the asymptotic distributions or simulations via fixing nuisance parameters or functions under the null hypothesis at certain values of interest. Furthermore, Fan, Zhang and Zhang (2001) showed that the resulting tests are asymptotically optimal in the sense of Ingster (1993).

The idea of the above generalized likelihood method is widely applicable in semiparametric and nonparametric models. It is easy to use because of the Wilks phenomenon and is powerful as it achieves the optimal rate of convergence. Yet, one needs to specify the parametric form of the error distribution such as $\varepsilon$ in (1.1) in order to construct the generalized likelihood ratio statistic. While the procedure based on the normal likelihood may be still applicable to the case where the distribution of $\varepsilon$ is homoscedastic, it may not be efficient. When the error distribution is heteroscedastic with the variance $\text{var}(\varepsilon|U) = \sigma^2(U)$, the construction of the generalized likelihood ratio test statistic needs the knowledge of the variance function $\sigma^2(\cdot)$. This motivates us to propose the sieve empirical likelihood ratio (SELR) test statistic for handling the case where the exact form of the error distribution is unknown, but some qualitative traits of the distribution are known. A popular model is to assume

$$(1.4) \qquad E[G(\varepsilon)|U] = 0$$

where $G = (G_1, \ldots, G_{k_0})^\tau$ is a $k_0$-dimensional function [see Owen (1990), Newey (1993) and Zhang and Gijbels (2003)]. This is a much less restrictive assumption than a parametric form on the distribution of $\varepsilon$. In particular, when the conditional distribution of $\varepsilon$ given $U$ is symmetric about 0, we may choose a sequence of $k_0$ grid points, say, $0 = s_0 < s_1 < \cdots < s_{k_0}$ and take

$$(1.5) \quad G_k(\varepsilon) = I(\varepsilon \in [s_{k-1}, s_k]) - I(-\varepsilon \in [s_{k-1}, s_k]), \qquad 1 \leq k \leq k_0,$$

or a smoother version of the function $G_k$, where $I(\cdot)$ is the indicator function. Note that as $\max_{1 \leq k \leq k_0}(s_k - s_{k-1}) \to 0$, $k_0 \to \infty$, these restrictions are essentially the same as the symmetric assumption on the distribution of $\varepsilon$.

A few questions related to the SELR test arise naturally. First of all, it is not clear how to construct an empirical likelihood in the nonparametric setting. Second, it is not obvious whether a particular construction of the empirical likelihood ratio statistic will follow the Wilks type of result. Third, it is not granted that the resulting test statistic is asymptotically optimal in the sense of Ingster (1993). Finally, it remains unknown whether the empirical likelihood ratio statistics will adapt to the unknown distribution of $\varepsilon$ including heteroscedasticity. These issues are poorly understood and need to be studied.

The technical derivations for SELR tests are very involved. To ease some of the technical burden, we choose the varying coefficient model (1.1) for our investigation. The model arises from various contexts and has been widely used. For example, in many biomedical studies one frequently encounters the issue of the extent to which the effect of exposure variables on the response variable changes with the level of a confounding covariate (e.g., age). See, for example, Cleveland, Grosse and Shyu (1991), Hastie and Tibshirani (1993) and Carroll, Ruppert and Welsh (1998). The model can also be used for predicting group behavior in economics where different groups are allowed to have different coefficients. In longitudinal studies, investigators often want to examine how the effects of covariates on response variables change over time [Brumback and Rice (1998) and Wu, Chiang and Hoover (1998)]. In nonlinear time series, the model allows different autoregressive models for different regimes of state variables [Chen and Tsay (1993) and Cai, Fan and Yao (2000)]. It includes the threshold autoregressive model [Tong (1990)] as a specific example. The model has successfully been applied by Hong and Lee (2003) to the inference and forecast of exchange rates. Thus, our study in model (1.1) has direct implications for the above problems.

For the varying coefficient model (1.1), whether the coefficient functions are really varying or whether certain covariates are statistically significant frequently arises. This leads to the problem of testing for homogeneity (1.2) or the problem of testing for significance such as the problem (1.3). As will be explained at the end of Section 2, these problems can be reduced to that of testing against a specific null hypothesis:

$$H_0 : a_1(\cdot) = a_{10}(\cdot), \ldots, a_p(\cdot) = a_{p0}(\cdot),$$

for some given functions $a_{10}, \ldots, a_{p0}$. Our approach is to first construct the local linear estimator of the coefficient functions $a_1, \ldots, a_p$ via a local version of the empirical likelihood, and to then substitute the estimate into a special sieve empirical likelihood [see Zhang and Gijbels (2003) and Zhang and Liu (2003)]. This allows us to form the empirical likelihood ratio statistics. We will show that the proposed SELR procedures follow the Wilks type of results under more relaxed assumptions on the error distribution of $\varepsilon$. This provides a useful extension of the results given by Fan, Zhang and Zhang (2001). Note that our procedure is very different from that of Kitamura (1997), who considered testing problems for finite-dimensional parameters in weakly dependent processes. He first used the local (blocking) approximation to construct a global estimating equation, then applied Owen's procedure directly. For the full nonparametric regression model, Chen, Härdle and Li (2003) developed a very different version of the empirical likelihood ratio test, using a kernel-smoothed parametric estimator under the null hypothesis as ancillary information. The idea has nicely been extended to simultaneously testing the parametric forms of the mean and variance functions by Chen, Gao and Li (2003). Horowitz and Spokoiny (2002) developed a different test for one special case of the model (1.4). The test shares most of nice features listed for the SELR test and includes an automatic selection of the smoothing parameter. It is not clear whether the Horowitz–Spokoiny test is adaptive to the error distribution under the alternative hypothesis, as is the SELR. Furthermore, because of the saturated alternative, the curse-of-dimensionality problem arises in implementation and power.

Our empirical likelihood ratio method applies also to the nonparametric tests on density functions. As an illustration without introducing new statistical setting, we regard the constraints (1.4) as a null hypothesis. We will demonstrate that the Wilks type of phenomenon continues to hold for this nonparametric testing problem.

Our studies have implications for other nonparametric models. When $p = 1$ and $X \equiv 1$, the model is the nonparametric regression model studied by many authors. Our results can be directly applied to the problems of testing parametric models against the nonparametric alternative. Further, our theoretical results shed some lights on the validity of the Wilks phenomenon in other models such as additive models and models under certain mixing conditions.

When $p = 1$ and $X \equiv 1$ and the coefficient function $a_1(\cdot) \equiv \theta$, under the constraints (1.4) and (1.5), the model becomes a one-sample symmetric location model, which was well studied, for instance, by Hettmansperger (1984) and Bickel, Klaassen, Ritov and Wellner (1993). In Section 2, we find that for this special case, the first step in our procedure essentially makes the information on the stochastic error to be efficiently used [Owen (1988) and Zhang and Liu (2003)]. Moreover, the second step makes the likelihood ratio statistic adaptive to heteroscedasticity. As a result, our procedure has two advantages over the parametric assumptions on the error distribution. First, it requires only some conditional estimating equations such as (1.4) rather than the whole distribution of the stochastic error. Second, the asymptotic null distribution of the SELR statistic asymptotically follows a rescaled $\chi^2$-distribution. The scaling constant and the degrees of freedom are independent of the conditional distribution of $\varepsilon$ even if the stochastic error is heteroscedastic in $U$. The procedure and results can be easily generalized to a more general constrained regression model in Zhang and Gijbels (2003).

The paper is organized as follows. In Section 2 the sieve empirical likelihood ratio statistics are introduced for testing the goodness-of-fit of the estimating equations and for testing some simple and composite null hypotheses. In Section 3 the asymptotic null and nonnull distributions of these statistics are derived. In Section 4 a simulation study is conducted to evaluate the performance of the proposed procedure empirically. The technical conditions and the proofs are relegated to Section 5. The technical lemmas are given in the Appendix.

**2. Sieve empirical likelihood.** It is more convenient to work with the matrix notation for model (1.1),

$$Y = A^\tau(U)X + \varepsilon, \tag{2.1}$$

where $Y$ is the response, $U \in \Omega \subset R^1$ (with $\Omega$ bounded) and $X \in R^p$ are covariates, $\varepsilon$ is the stochastic error and $A(U) = (a_1(u), \ldots, a_p(u))$ is the vector of varying coefficients. Let $\{(Y_i, X_i, U_i)\}_{i=1}^n$ be an i.i.d. random sample from the model (2.1) with the restriction (1.4). According to Owen (1990), to construct an empirical likelihood which can identify an infinite-dimensional parameter such as $A(u)$ in (2.1), we need to establish an infinite number of unconditional estimating equations. Such a likelihood is often theoretically intractable. To overcome this difficulty, Zhang and Gijbels (2003) proposed a general procedure to build a sieve empirical likelihood via local approximation. For the model (2.1) the procedure consists of two steps: First, for each $U_j$ construct $n$ local empirical likelihoods which can locally identify $A(u), u \approx U_j$. These local empirical likelihoods lead to a weighted approximation of the logarithm of the conditional probability mass $dP_{(Y,X)|U=U_j}(Y_j, X_j)$. Then a log-likelihood is obtained simply by summing up all of these approximated logarithms. In the first step, we will use the local linear approximation of the nonparametric coefficient functions $A(\cdot)$ [see Fan (1992), Fan and Zhang (1999) and Cai, Fan and Li (2000)]. In other words, in a neighborhood around a given point $u_0$, approximate $A(\cdot)$ by

$$A(u) \approx A(u_0) + A'(u_0)(u - u_0) \quad \text{for } u \approx u_0.$$

Thus, around the point $u_0$, the model (2.1) and the restriction (1.4) can be written, respectively, as

$$
\begin{aligned}
Y &\approx \beta_A(u_0)^T Z\left(X, \frac{U - u_0}{h}\right) + \varepsilon \quad \text{for } U \approx u_0, \\
E\left[G(Y - \beta_A(u_0))^\tau Z\left(X, \frac{U - u_0}{h}\right) \Big| U = u\right] &\approx 0 \quad \text{for } u \approx u_0,
\end{aligned}
\tag{2.2}
$$

where $\beta_A(u_0) = (A^\tau(u_0), hA'^\tau(u_0))^\tau$ and $Z(X, t) = (X^\tau, tX^\tau)^\tau$. This is indeed a local linear model [Fan (1992)]. To incorporate the local linear model, let $h$ represent the size of the local neighborhood where the approximation is valid and $K$ be a weight function, which is a symmetric probability density function. Let $p_i, i = 1, \ldots, n$, be the conditional empirical probability mass of $(X, Y)$ given $U = u_0$, putting on the $i$th data point $(X_i, Y_i)$, $i = 1, \ldots, n$.

Suppose that given $U$, $\varepsilon$ and $X$ are independent. Then the conditional constraints (2.2) can be translated into the following unconditional estimating equation:

$$\sum_{i=1}^{n} p_i \mathbf{G}_{ih}(u_0, \beta_A(u_0)) = 0,$$

where

$$\mathbf{G}_{ih} = \mathbf{G}_{ih}(u_0, \beta) = G\bigg(Y_i - \beta^\tau Z\bigg(X_i, \frac{U_i - u_0}{h}\bigg)\bigg) \otimes Z\bigg(X_i, \frac{U_i - u_0}{h}\bigg)$$

with $\otimes$ being the Kronecker product, $\beta = (A^\tau, hB^\tau)^\tau$, $A = (a_1, \ldots, a_p)^\tau$ and $B = (b_1, \ldots, b_p)^\tau$. To see why we need an extra factor $Z(X_i, (U_i - u_0)/h)$ in the unconditional estimating function $\mathbf{G}_{ih}$, we let $G(\varepsilon) = \varepsilon$ temporarily. It is a well-known fact that in the linear model the product of the residual and the covariates is a good estimating equation for the parameter $\beta_A$. This leads to the estimating equation:

$$\sum_{i=1}^{n} p_i \bigg(Y_i - \beta^\tau Z\bigg(X_i, \frac{U_i - u_0}{h}\bigg)\bigg) Z\bigg(X_i, \frac{U_i - u_0}{h}\bigg) = 0.$$

In light of this fact, for a general $G$ we should build the estimating equation by multiplying each component of $G$ by the covariate vector $Z(X_i, (U_i - u_0)/h)$, which admits the form $\mathbf{G}_{ih}$.

Thus, following Owen (1988, 1990), the local empirical log-likelihood function of $\beta$ is defined by

(2.3)
$$l(\beta, u_0) = \sup\bigg\{\sum_{i=1}^{n} w_h(U_i, u_0) \log p_i : p_i \geq 0, 1 \leq i \leq n,$$
$$\sum_{i=1}^{n} p_i = 1, \sum_{i=1}^{n} p_i \mathbf{G}_{ih}(u_0, \beta) = 0\bigg\},$$

where $w_h(U_i, u_0) = K_h(U_i - u_0)/\sum_{m=1}^{n} K_h(U_m - u_0)$ with $K_h(\cdot) = K(\cdot/h)/h$. If we set $p_i = w_h(U_i, u_0) q_i$, then (2.3) becomes

$$l(\beta, u_0) = \sup\bigg\{\sum_{i=1}^{n} w_h(U_i, u_0) \log\{w_h(U_i, u_0) q_i\} : q_i \geq 0, 1 \leq i \leq n,$$
$$\sum_{i=1}^{n} w_h(U_i, u_0) q_i = 1, \sum_{i=1}^{n} w_h(U_i, u_0) q_i \mathbf{G}_{ih}(u_0, \beta) = 0\bigg\}.$$

Analogously to Owen (1990) and Qin and Lawless (1994), if 0 is contained in the convex hull of the points in $\{\mathbf{G}_{ih}(u_0, \beta) : w_h(U_i, u_0) > 0, 1 \leq i \leq n\}$, then an explicit expression can be derived by the Lagrange multiplier method as follows:

$$l(\beta, u_0) = \sum_{i=1}^{n} w_h(U_i, u_0) \log w_h(U_i, u_0)$$
$$- \sum_{i=1}^{n} w_h(U_i, u_0) \log(1 + \alpha_n^\tau(u_0, \beta)\mathbf{G}_{ih}(u_0, \beta)),$$

where $\alpha_n(u_0, \beta)$ satisfies

$$(2.4) \qquad \sum_{i=1}^{n} w_h(U_i, u_0) \frac{\mathbf{G}_{ih}(u_0, \beta)}{1 + \alpha_n^\tau(u_0, \beta)\mathbf{G}_{ih}(u_0, \beta)} = 0.$$

Define the estimate of $\beta$ by

$$(2.5) \qquad \hat{\beta}(u_0) = \arg\max_{\beta} l(\beta, u_0).$$

The first $p$ components, denoted by $\hat{A}(u_0)$, give an estimate of $A(u_0)$, and the remaining components estimate $hA'(u_0)$. Similarly to LeBlanc and Crowley (1995), an approximate empirical likelihood, called the sieve empirical likelihood for the nonparametric function $A$ can be introduced by adding the logarithm of the conditional likelihood at each data point:

$$l(A|G) = \sum_{j=1}^{n} l(\beta_A, U_j).$$

The name "sieve" originates from the following two facts: $\{E[G(\varepsilon)|U = U_j]\}_{1 \leq j \leq n}$ is a sieve approximation to the constraints (1.4) and $l(\beta_A, U_j)$ is a weighted approximation of the logarithm of the conditional probability mass $dP_{(Y,X)|U=U_j}(Y_j, X_j)$. See Zhang and Gijbels (2003) for a more detailed explanation. Motivated by Fan, Zhang and Zhang (2001), we define the logarithm of the sieve empirical likelihood under the nonparametric model (2.1) with constraints (1.4) by substituting $\beta = \hat{\beta}$ into $l(A|G)$, leading to

$$l(\Theta|G) = \sum_{j=1}^{n} l\{\hat{\beta}(U_j), U_j\}$$

with $\Theta$ denoting the space of $A$.

We now consider the nonparametric test concerning the density function of $\varepsilon$. As a specific application of the sieve empirical likelihood, we consider testing

$$(2.6) \qquad H_{0G} : E[G(\varepsilon)|U] = 0,$$

where $G$ is given in (1.4). Without the constraint (1.4), following the above derivations, the corresponding logarithm of the sieve empirical likelihood is

$$l(\Theta|N) = \sum_{j=1}^{n} \sum_{i=1}^{n} w_h(U_i, U_j) \log w_h(U_i, U_j).$$

Thus, we can construct a goodness-of-fit test of hypothesis (2.6) based on the following logarithm of the SELR:

$$\begin{aligned}
l(G) &= -l(\Theta|G) + l(\Theta|N), \\
&= \sum_{j=1}^{n} \sum_{i=1}^{n} w_h(U_i, U_j) \log(1 + \hat{\alpha}(U_j)^\tau \mathbf{G}_{ih}(U_j, \hat{\beta}))
\end{aligned} \quad (2.7)$$

where $\hat{\alpha}(u) = \alpha_n(u, \hat{\beta})$.

Next, we consider the sieve likelihood ratio test for the nonparametric coefficient function $A(\cdot)$ under the restriction (1.4). In the varying coefficient model (2.1), we ask naturally whether the coefficient is really varying or whether certain covariates are statistically significant. This leads to the parametric null hypothesis:

$$H_{0p}: A(\cdot) = \theta.$$

More generally, we wish to test the composite null hypothesis, which involves nuisance functions $A_2(\cdot)$:

$$(2.8) \qquad H_{0u}: A_1 = A_{10} \quad \Longleftrightarrow \quad H_{1u}: A_1 \neq A_{10}$$

with $A_2(\cdot)$ completely unknown. This problem includes the test of significance (1.3) under model (1.1) as a specific example. Here we write

$$A_0(u_0) = \begin{pmatrix} A_{10}(u_0) \\ A_{20}(u_0) \end{pmatrix} \quad \text{and} \quad A(u_0) = \begin{pmatrix} A_1(u_0) \\ A_2(u_0) \end{pmatrix},$$

with $A_{10}(u)$ and $A_1(u)$ being $p_1$ ($< p$)-dimensional. To construct the likelihood ratio statistic for $H_{0u}$, we introduce the following notation:

$$\beta_{2A}(u_0) = (A_2^\tau(u_0), hA_2'^\tau(u_0)), \qquad \beta_2 = (A_2^\tau, hB_2^\tau)^\tau,$$
$$\beta^* = (A_{10}^\tau(u_0), A_2^\tau, hA_{10}'^\tau(u_0), hB_2^\tau)^\tau.$$

Let

$$\hat{\beta}_2(u_0) = (\hat{A}_2^\tau, h\hat{B}_2^\tau)^\tau = \arg\max_{\beta_2} l(\beta^*, u_0),$$

$$\hat{\beta}^*(u_0) = (A_{10}^\tau(u_0), \hat{A}_2^\tau, hA_{10}'^\tau(u_0), h\hat{B}_2^\tau)^\tau$$

and the corresponding $\hat{\alpha}^*(u_0)$ be implicitly defined by

$$0 = \frac{1}{n} \sum_{i=1}^{n} w_h(U_i, u_0) \frac{\mathbf{G}_{ih}(u_0, \hat{\beta}^*(u_0))}{1 + \hat{\alpha}^{*\tau}(u_0) \mathbf{G}_{ih}(u_0, \hat{\beta}^*(u_0))}.$$

Then the SELR statistic for $H_{0u}$ can be written as

$$l(H_{0u}|G) = -l(\Theta_{02}|G) + l(\Theta|G), \tag{2.9}$$

where $\Theta_{02}$ denotes the space of $A_2$ and

$$l(\Theta_{02}|G) = \sum_{j=1}^{n} l(\hat{\beta}^*(U_j), U_j).$$

The SELR test for the semiparametric model that $A(\cdot)$ has a certain parametric form such as the linear model can be constructed analogously. As in Fan, Zhang and Zhang (2001), the asymptotic null distributions of the SELR statistics for composite null hypotheses can be derived from those for simple hypotheses (see the next paragraph). This motivates us to consider

$$H_{0s}: A = A_0 \quad \Longleftrightarrow \quad H_{1s}: A \neq A_0 \tag{2.10}$$

for a given $A_0$. Analogously to $l(H_{0u}|G)$, we can construct the following likelihood ratio statistic:

$$\begin{aligned}
l(H_{0s}|G) &= -l(A_0|G) + l(\Theta|G) \\
&= \sum_{j=1}^{n}\sum_{i=1}^{n} w_h(U_i, U_j) \log(1 + \alpha_n(U_j, \beta_0)^\tau \mathbf{G}_{ih}(U_j, \beta_0)) - l(G)
\end{aligned} \tag{2.11}$$

where $\beta_0$ denotes $\beta_{A_0}$. Note that when $A_0$ in $H_{0s}$ is known, we can assume, without loss of generality, that $A_0 \equiv 0$. This can be accomplished through a simple transformation $A^* = A - A_0$. With this transformation, (2.10) is equivalent to

$$H'_{0s}: A^* \equiv 0 \quad \Longleftrightarrow \quad H'_{1s}: A^* \not\equiv 0. \tag{2.12}$$

This specific problem has an advantage: the local linear estimator under the null hypothesis is unbiased and hence the null distribution can be more accurately approximated.

We opt for general $A_0$, since the results have implications for the composite null hypotheses. To appreciate this, consider the composite null hypothesis testing problem:

$$H_0: A \in \mathcal{A}_0 \quad \Longleftrightarrow \quad A \notin \mathcal{A}_0, \tag{2.13}$$

where $\mathcal{A}_0$ is a set of functions. Let $l(\mathcal{A}_0|G)$ be the sieve empirical likelihood under the hypothesis $H_0$ in (2.13). Then, the SELR statistic is simply

$$\lambda_n = -l(\mathcal{A}_0|G) + l(\Theta|G).$$

Let $A'_0$ denote the true value of the parameter function $A$. Consider the fabricated testing problems with the simple null hypotheses:

$$H'_0: A = A'_0 \quad \Longleftrightarrow \quad H_1: A \neq A'_0 \tag{2.14}$$

and

(2.15) $$H_0': A = A_0' \iff H_1': A \in \mathcal{A}_0.$$

Let $l(A_0'|G)$ be the sieve empirical likelihood under $H_0'$. Then the SELR statistic for (2.13) can be written as

$$\lambda_n = \lambda(A_0'|G) - \lambda^*(A_0'|G),$$

where $\lambda(A_0'|G) = -l(A_0'|G) + l(\Theta|G)$ is the SELR statistic for the problem (2.14) and $\lambda^*(A_0'|G) = -l(A_0'|G) + l(\mathcal{A}_0|G)$ is the SELR test for the problem (2.15). Thus, the asymptotic representation of $\lambda_n$ follows directly from those of $\lambda(A_0')$ and $\lambda^*(A_0')$, which admits the form given by (2.11).

## 3. Asymptotic theory.

3.1. *Asymptotic expansions.* In order to obtain the properties of the SELR statistics in (2.7) and (2.11), we first develop some uniform asymptotic representations for the local sieve empirical likelihood estimator $\hat{\beta}$ and the Lagrange multiplier $\hat{\alpha}$ in (2.4) and (2.5). These results are the generalizations of Zhang and Liu (2003). They also indicate the performance of the sieve empirical likelihood estimator. Using these results we will establish the asymptotic representations for $l(G)$ and $l(H_{0s}|G)$ in (2.7) and (2.11). For simplicity of presentation, we assume $G$ is differentiable. Let $f(u_0)$ be the density of $U$ at the point $u_0$. Set

$$D(u_0) = -E\left[\frac{\partial G(\varepsilon)}{\partial \varepsilon}\Big| U = u_0\right],$$

$$V(u_0) = E[G(\varepsilon)G^\tau(\varepsilon)|U = u_0],$$

$$\Gamma(u_0) = E[XX^\tau|U = u_0]f(u_0),$$

$$S = \begin{pmatrix} 1 & 0 \\ 0 & \mu_2 \end{pmatrix}, \qquad \mu_2 = \int t^2 K(t)\,dt,$$

$$\eta_i(u_0) = -\{D(u_0)^\tau V(u_0)^{-1} D(u_0)\}^{-1} D^\tau(u_0) V^{-1}(u_0) G(\varepsilon_i),$$

$$C(u_0) = V^{-1}(u_0) - V^{-1}(u_0)(D^\tau(u_0)V^{-1}(u_0)D(u_0))^{-1}D(u_0)D^\tau(u_0)V^{-1}(u_0),$$

$$\varepsilon_i = Y_i - A^\tau(U_i)X_i.$$

THEOREM 1. *Suppose that conditions* (K0), (U0), (A1)–(A10) *and* (B1)–(B5) *in Section* 5.1 *hold and that the underlying* $A(u)$ *have twice continuous derivatives and satisfy condition* (B6). *If there exist some positive constants* $b_0, b_1$ *and* $\eta < 1/2$ *such that* $b_0 \leq hn^\eta \leq b_1$, *then uniformly for* $u_0 \in \Omega$,

$$\hat{\beta}(u_0) = \beta(u_0) + \frac{1}{n}\sum_{i=1}^n K_h(U_i - u_0) \begin{pmatrix} \Gamma^{-1}(u_0)X_i \\ \mu_2^{-1}\Gamma^{-1}(u_0)X_i(U_i - u_0)/h \end{pmatrix}$$

$$\times \eta_i(u_0)(1 + o_p(h^{1/2})) + O_p(h^2),$$

$$\hat{\alpha}(u_0) = \frac{1}{n}\sum_{i=1}^{n} K_h(U_i - u_0)\{C(u_0)G(\varepsilon_i)\}$$

$$\otimes \begin{pmatrix} \Gamma^{-1}(u_0)X_i \\ \mu_2^{-1}\Gamma^{-1}(u_0)X_i(U_i - u_0)/h \end{pmatrix}(1 + o_p(h^{1/2})) + O_p(h^2).$$

As a consequence of Theorem 1, we have the following asymptotic uniform expansion:

$$\hat{A}(u_0) - A(u_0) = \frac{1}{n}\sum_{i=1}^{n} K_h(U_i - u_0)\Gamma^{-1}(u_0)X_i\eta_i(u_0)(1 + o_p(h^{1/2})) + O_p(h^2).$$

The asymptotic normality of the local sieve empirical likelihood estimator follows easily from the above asymptotic expansion.

In Theorem 1, the requirement that $G$ is differentiable can be relaxed by imposing some entropy conditions on $G$ and by assuming $E[G(\varepsilon - t)|U = u_0]$ is twice continuously differentiable in $t$. In this case $D(u_0)$ should be replaced by $-\{\partial E[G(\varepsilon - t)|U = u_0]/\partial t\}|_{t=0}$. Similarly to Zhang and Liu (2003), we can show that the asymptotic efficiency of $\hat{A}(u_0)$ is increasing in $D(u_0)^\tau V(u_0)^{-1} \times D(u_0)$. In particular, in the setting of the symmetric location model mentioned in Section 1, we can find a sequence of $G$ functions, say $\{G^{(k)}\}$, such that the corresponding $\hat{A}(u_0)$ is asymptotically adaptive to the unknown conditional density of $\varepsilon$ given $U = u_0$. In practice, to save computational effort, we prefer to choose a $G$ with a small $k_0$ and a relatively large $D(u_0)V(u_0)^{-1}D(u_0)$.

It should be noted that under the conditions of Theorem 1, for $\beta$ near its true value, $\alpha(u_0, \beta)$ is uniquely determined by the estimating equations. Thus, the number of unknown parameters is $2p$ for each $u_0$. It is well known that to make the local linear model regular, the interval $[u_0 - h, u_0 + h]$ should include at least $2p + 1$ data points of $U$. This condition asymptotically holds under the condition of Theorem 1 because as $n \to \infty$,

$$P\left\{\sum_{i=1}^{n} I(u_0 - h \leq U_i \leq u_0 + h) \geq 2p + 1\right\}$$

$$= P\left\{\left|\sum_{i=1}^{n}[I(u_0 - h \leq U_i \leq u_0 + h)\right.\right.$$

$$- EI(u_0 - h \leq U \leq u_0 + h)]$$

$$\left.\left. + nEI(u_0 - h \leq U \leq u_0 + h)\right| \geq 2p + 1\right\}$$

$$\geq P\{nEI(u_0 - h \leq U \leq u_0 + h) \geq 2p + 1 + \delta\}$$

$$-nE(I(u_0 - h \leq U \leq u_0 + h) - EI(u_0 - h \leq U \leq u_0 + h))^2/\delta^2$$
$$\to 1$$

where $\delta = n^{(1+2\eta)/2}h$ and $I(\cdot)$ is the indicator function. We can further show that this condition actually holds uniformly in $u_0$ by an approach using empirical processes.

We now give the asymptotic representations for the SELR statistics $l(G)$ and $l(H_{0s}|G)$. The results indicate that they admit a generalized quadratic form. To facilitate the expressions, the following notation is introduced. Let

$$\phi_{ikh}(U) = K_h(U_i - U)K_h(U_k - U)C(U)$$
$$\times (1 + (U_i - U)(U_k - U)\mu_2^{-1}h^{-2})X_i^\tau \Gamma^{-1}(U)X_k f^{-1}(U),$$

(3.1) $$K^*(s) = \int K(t)K(s+t)(1 + t(s+t)\mu_2^{-1})\,dt,$$

$$\Phi_{ikh} = E[\phi_{ikh}(U)|(U_i, U_k, X_i, X_k)]$$
(3.2) $$= K_h^*(U_k - U_i)C(U_i)X_i^\tau \Gamma^{-1}(U_i)X_k(1 + O_p(h)),$$

$$T_n = \frac{1}{n(n-1)} \sum_{i \neq k} G^\tau(\varepsilon_i)\Phi_{ikh}G(\varepsilon_k).$$

Similarly, we define

$$q_{ikh}(U) = K_h(U_i - U)K_h(U_k - U)V^{-1}(U)X_i^\tau \Gamma^{-1}(U)$$
$$\times X_k\{1 + (U_i - U)(U_k - U)\mu_2^{-1}h^{-2}\}f^{-1}(U),$$

$$Q_{ikh} = E[q_{ikh}(U)|(U_i, U_k, X_i, X_k)],$$

$$T_n^* = \frac{1}{n(n-1)} \sum_{i \neq k} G^\tau(\varepsilon_i)(Q_{ikh} - \Phi_{ikh})G(\varepsilon_k).$$

Then we have the following result.

THEOREM 2. *Suppose the conditions of Theorem* 1 *hold. Then under* $H_{0G}$,

(3.3) $$2l(G) = \frac{(k_0 - 1)p|\Omega|}{h} \int K^2(t)(1 + t^2\mu_2^{-1})\,dt + (1 + o_p(h^{1/2}))nT_n + o_p(h^{-1/2});$$

*and under* $H_{0s}$, *if* $A_0$ *is linear or* $nh^{9/2} \to 0$, *then*

(3.4) $$2l(H_{0s}|G) = \frac{p|\Omega|}{h} \int K^2(t)(1 + t^2\mu_2^{-1})\,dt + (1 + o_p(h^{1/2}))nT_n^* + o_p(h^{-1/2}),$$

*where* $|\Omega|$ *is the length of the support* $\Omega$ *of the density* $f$.

Note that if there are no components in $A$, then under $H_{0G}$ the factor $k_0 - 1$ in (3.3) should be $k_0$, since it costs $p$ degrees of freedom to estimate them when there are $p$ components in $A$.

3.2. *Asymptotic null distribution.* With the asymptotic representations, we are now ready to derive the asymptotic distributions of the test statistics $l_G$ and $l(H_{0s}|G)$. As in the parametric case for the stochastic error $\varepsilon$ [see Fan, Zhang and Zhang (2001)], under the null hypotheses the SELR statistics in (2.7), (2.9) and (2.11) are asymptotically $\chi^2$-distributed and their degrees of freedom are independent of the nuisance parameters such as $A$, $G$ and the distribution of $\varepsilon$.

THEOREM 3. *Under $H_{0G}$ and the conditions of Theorem 1, for $k_0 > 1$, we have $r_K l_G \stackrel{a}{\sim} \chi^2_{b_n}$ with*

$$r_K = \frac{2K^*(0)}{\int K^*(s)^2 \, ds}, \qquad b_n = \frac{(k_0 - 1)p|\Omega|c_K}{h},$$

*where $K^*(s)$ is defined in (3.1), $c_K = K^*(0)^2 / \int K^*(s)^2 \, ds$. For $k_0 = 1$, we have $r_K l_G = o_p(1)$.*

REMARK 1. If $K(t)$ has support $[-1, 1]$, and if $K(t)$ and $|t|K(t)$ are concave on $t \in [-1, 1]$, then by the same argument used in the Sherman inequality [see Farrell (1985), page 343], we have

$$|K^*(s)| \leq \int K(t)K(s+t) \, dt + \mu_2^{-1} \int |t|K(t)|s+t|K(s+t) \, dt$$
$$\leq K^*(0).$$

Thus when $K^*(s) \geq 0, s \in [-1, 1]$, $r_K \geq 2$. In particular, when $K$ is the uniform kernel function, $r_K = 2.8176$ and $c_K = 1.0566$; when $K$ is the Epanechnikov kernel function, $r_K = 2.5154$ and $c_K = 1.2936$.

The next theorem presents the asymptotic null distribution of $l(H_{0s}|G)$.

THEOREM 4. *Suppose that the conditions of Theorem 1 hold. Then under $H_{0s}$, $r_K l(H_{0s}|G) \stackrel{a}{\sim} \chi^2_{b_n^*}$ if $A_0$ is linear or $nh^{9/2} \to 0$; and under $H_{0u}$, if $nh^{9/2} \to 0$, then $r_K l(H_{0u}|G) \stackrel{a}{\sim} \chi^2_{b_{n2}^*}$ where $b_n^* = p|\Omega|c_K/h$ and $b_{n2}^* = p_1|\Omega|c_K/h$ with $c_K$ and $r_K$ defined in Theorem 3 and $p_1$ being the dimensionality of $A_{10}$ in (2.8).*

When $nh^{9/2} = O(1)$, it is easily proved as in Fan, Zhang and Zhang (2001) that under $H_{0u}$ the Wilks phenomenon continues to hold in the generalized sense that the mean and variance of the SELR statistic are independent of the nuisance parameters to the first order. As pointed out in Section 2, when $A_0$ in $H_{0s}$ is known (or more generally in a parametric form), we can make a simple transformation (or use some bias reduction technique) to kill the bias. Theorems 3 and 4 indicate that the SELR statistics continue to apply to the case where the distribution of the stochastic error $\varepsilon$ is completely unknown and, furthermore, there are many nuisance parameters in null hypotheses (see Section 3.4). In particular, the stochastic errors are allowed to be heteroscedastic and unknown. This is a useful generalization of the results in Fan, Zhang and Zhang (2001) where the distribution of $\varepsilon$ is essentially known. In particular, if the variance is heteroscedastic with $\mathrm{var}(\varepsilon|U) = \sigma^2(U)$, they have to rely on the knowledge of $\sigma^2(\cdot)$ to construct the likelihood ratio statistics. This drawback is repaired by the empirical likelihood ratio method, while their Wilks phenomenon is inherited.

3.3. *Asymptotic power.* To demonstrate the effectiveness of the sieve empirical likelihood method, we consider, for simplicity, the test statistic for the problem (2.12) under the contiguous alternative $A_n(\cdot) \to 0$, with $A_n''(\cdot)$ being bounded; that is, we allow the coefficient functions to be close to the null hypothesis, but still in the class of functions with bounded and continuous second derivatives. This is a much weaker restriction than the contiguous alternatives of the form $A_n(u) = a_n B_0(u)$ for a sequence $a_n \to 0$ and a given $B_0$, considered by many authors [e.g., Eubank and Hart (1992), Eubank and LaRiccia (1992), Hart (1997) and Inglot and Ledwina (1996)]. The latter implicitly assumes that $A_n'(u) \to 0$ and $A_n''(u) \to 0$, which are too restrictive for nonparametric applications.

We begin with the following notation. Let

$$W_{1n}^* = \frac{1}{n} \sum_{i \neq k}^n K_h^*(U_i - U_k) G(\varepsilon_i) V^{-1}(U_k) \tag{3.5}$$
$$\times X_i^\tau \Gamma^{-1}(U_k) X_k A(U_k)^\tau X_k \frac{\partial G(\varepsilon_k)}{\partial \varepsilon},$$

$$\Xi_i = \frac{\partial G(\varepsilon_i)}{\partial \varepsilon} - E\left[\frac{\partial G(\varepsilon_i)}{\partial \varepsilon}\Big|U_i\right], \tag{3.6}$$

$$W_{2n}^* = \frac{1}{n} \sum_{i \neq k}^n K_h^*(U_i - U_k) \Xi_i^\tau V^{-1}(U_i) \Xi_k A(U_i)^\tau \tag{3.7}$$
$$\times X_i X_i^\tau \Gamma^{-1}(U_k) X_k X_k^\tau A(U_k),$$

$$W_{3n}^* = \frac{1}{n} \sum_{i \neq k}^n K_h^*(U_i - U_k) \Xi_i^\tau V^{-1}(U_k) E\left[\frac{\partial G(\varepsilon_k)}{\partial \varepsilon}\Big|U_k\right] \tag{3.8}$$
$$\times A(U_i) X_i X_i^\tau \Gamma^{-1}(U_k) X_k X_k^\tau A(U_k).$$

Then, following the same arguments used in Fan, Zhang and Zhang (2001), we can derive the asymptotic power $l(H_{0s}|G)$ via the next theorem.

THEOREM 5. *Assume that $A_0 \equiv 0$ and that the underlying coefficient $A = A_n$ has twice continuous derivatives and satisfies $nhEA(U)^\tau XX^\tau A(U) = O(1)$, $\max_u \|A(u)\| \to 0$ and $\max_u \|A''(u)\| = O(1)$ as $n \to \infty$. Assume that $G$ is twice continuously differentiable. Then under the conditions of Theorem 1,*

$$2l(H_{0s}|G)$$
$$= \frac{p|\Omega|}{h} K^*(0) + nE\{D(U)^\tau V^{-1}(U)D(U)A(U)^\tau XX^\tau A(U)\}(1+o(1))$$
$$- \frac{nh^4}{4} E\{D(U)^\tau C(U)D(U)A''(U)^\tau XX^\tau A''(U)\}$$
$$\times \int\int t^2(s+t)^2 K(t)K(s+t)(1+\mu_2^{-1}t(s+t))\,dt\,ds\,(1+o(1))$$
$$+ (1+o_p(h^{1/2}))\{T_n^* + 2W_{1n}^* + W_{2n}^* + 2W_{3n}^*\} + o_p(h^{-1/2}),$$

where $D$, $V$, $C$ and $K^*$ are defined in Section 3.1.

Using the above result, similar to that in Fan, Zhang and Zhang (2001), it can easily be shown that under $H_{0s}$ the SELR can detect the alternative with rate $n^{-4/9}$ when $h = c_* n^{-2/9}$ for some constant $c_*$. This rate is optimal in the ordinary nonparametric regression setting. Note that the above result continues to hold for the composite null hypothesis testing problem (2.13) when $\mathcal{A}_0$ is a set of linear functions.

3.4. *Remarks on practical implementations.* There are a couple of issues arising from practical implementations of the procedure, including computing $P$-values, choice of bandwidths, choice of the support of $U$ and bias reduction. We now briefly discuss them.

$P$-values depend on the null distributions of test statistics. The convergence of the null distributions of the SELR statistics is expected to be slow. Thus, we do not suggest using the asymptotic null distributions. Instead, we use simulation methods (a form of bootstrap). Thanks to Theorems 3 and 4, we can simulate the null distributions by fixing nuisance parameters or functions under the null hypothesis at certain values of interest. This will give better approximations to the null distributions. We have conducted an intensive simulation study in Section 4. The results show that for a sample size of 200 or more, the approach gives very reasonable approximations of the null distribution.

The SELR test depends on the choice of bandwidth $h$. It can be regarded as a family of test statistics indexed by the bandwidth $h$. A thorough discussion of this subject is beyond the scope of this study. Inspired by the adaptive Neyman test in Fan (1996), which has been demonstrated to be adaptive minimax by Fan and Huang (2001), one can possibly use the following criterion to choose a bandwidth: For some constants $a, b > 0$, a bandwidth $\hat{h} \in [n^{-a}, n^{-b}]$ is selected to give a maximum value of

$$\frac{r_0 l(H_{0s}|G) - d_n(h)}{\sqrt{2 d_n(h)}},$$

where $r_0$ is the normalizing constant and $d_n(h)$ is the degrees of freedom (see Theorem 4). This results in a multi-scale test:

$$\frac{r_0 l(H_{0s}|G) - d_n(\hat{h})}{\sqrt{2 d_n(\hat{h})}} = \max_{h \in [n^{-a}, n^{-b}]} \frac{r_0 l(H_{0s}|G) - d_n(h)}{\sqrt{2 d_n(h)}}.$$

Such an idea was proposed in Fan, Zhang and Zhang [(2001), page 175] and in Horowitz and Spokoiny (2002) for the median regression problem and was shown to possess the adaptive optimal rate of convergence [Horowitz and Spokoiny (2002)]. It has also been studied and implemented by Zhang (2003). In many empirical applications, the bandwidths used for nonparametric function estimation have also been frequently employed for nonparametric hypothesis testing. The difference between the optimal bandwidth $O(n^{-1/5})$ for function estimation and $O(n^{-2/9})$ for hypothesis testing is hardly noticeable for practical sample sizes.

When $U$ has an unbounded support, we can not estimate the coefficient functions $A(\cdot)$ at the tails with reasonably good accuracy. In other words, we do not have enough data to test on the form of the coefficient functions at the tails. Due to this limitation, a reduced problem needs to be considered: test on the form of the coefficient functions $A(\cdot)$ on a given interval. Our procedures continue to apply and $|\Omega|$ becomes the length of the given interval.

When $A_0$ in (2.10) is of parametric form $A(\cdot, \theta)$ and is nonlinear, the local linear estimate will be biased even under the null hypothesis. The bias is killed by requiring $nh^{9/2} \to 0$ in the second part of Theorem 4. This is an unrealistic assumption, as pointed out by a referee. However, as discussed in Section 2, we should employ a bias reduction technique before applying the SELR test. Let $\hat{\theta}$ be a root-$n$ consistent estimator under the parametric model. The error of parametric fit is usually negligible in nonparametric applications. By regarding $A(\cdot, \hat{\theta})$ as $A_0$ in (2.10), we can deduce the problem to (2.12). For problem (2.12), the local linear fit does not have any bias under the null hypothesis. Hence, the condition $nh^{-9/2} \to 0$ is not required to kill bias. The bias reduction is also helpful in reducing approximation errors of the null distribution.

In summary, for practical implementations, the following steps are recommended:

1. Apply the bias correction method as in the last paragraph.
2. Choose an interval where functions are to be tested. This is the set $\Omega$.
3. Choose an appropriate bandwidth, using the methods suggested above to construct a SELR.
4. Apply the bootstrap method above to obtain a null distribution of the test statistic.

**4. Simulation.** In this section the performance of the SELR test is evaluated for a simplified conditional regression model by simulation. In this study, several bandwidths (i.e., $h = c_0 n^{-2/9}$, with $c_0 = 1$ and 1.5 for the sample size $n = 100$; with $c_0 = 0.5, 1, 1.5$ and 2 for $n = 200$ and 400; with $c_0 = 0.55, 1, 1.5$ and 2 for $n = 800$; and with $c_0 = 0.2, 0.35, 0.55$ and 2 for $n = 1600$) are used to represent widely varying degrees of smoothness. Due to space limitation, only part of the results is presented. The triweight function $(1-t^2)^3_+$ is selected as the kernel function in the proposed test.

For simplicity of exposition and computation, we take the simple model,

$$Y = a_1(U) + \varepsilon,$$

where $E[\varepsilon|U] = 0$ [i.e., $G = \varepsilon$ in (1.4)] and $U$ is uniformly distributed over $[0, 1]$, though the results hold for more general varying-coefficient models. Consider the problem of nonparametric testing of significance:

$$H_0 : a_1(\cdot) \equiv 0.$$

The SELR test of $H_0$ can be expressed as

$$l(\Theta_{0s}|G) = \sum_{i=1}^n \sum_{j=1}^n w_h(U_i, U_j) \log(1 + \alpha(U_j, 0)^\tau \mathbf{G}_{ih}(U_j, 0)),$$

where $\alpha(U_j, 0)$ satisfies

$$\sum_{i=1}^n w_h(U_i, U_j) \frac{\mathbf{G}_{ih}(U_j, 0)}{1 + \alpha_n^\tau(U_j, 0)\mathbf{G}_{ih}(U_j, 0)} = 0$$

with

$$\mathbf{G}_{ih}(U_j, 0) = Y_j \times \left(1, \frac{U_i - U_j}{h}\right)^\tau.$$

Note that Theorems 3 and 4 imply that the null distribution of $l(\Theta_{0s}|G)$ is asymptotically independent of the underlying distribution of $\varepsilon$. So without loss of generality, we assume that given $U$ the stochastic error follows a normal distribution.

To examine the effect of the possible heteroscedasticity of $\varepsilon$ on the above null distributions, the conditional variance of $\varepsilon$ is taken to have the form $(1 + c_1 U^2)$, where the constant $c_1$ represents the noise level. By generating 100 independent samples of $(Y, U)$ with sample size $n$, we calculate the null distributions for several quite different values of $c_1$, which represent widely varying degrees of heteroscedasticity of $\varepsilon$. This results in 100 i.i.d. simulated values of $l(\Theta_{0s}|G)$ for each combination of $n$ and $c_1$. The corresponding sample means and variances of $l(\Theta_{0s}|G)$ summarize the distributions of the test statistics under the null hypothesis and are reported in Table 1. They do not strongly depend on the choice of the constant $c_1$. As an illustration, the resulting 24 empirical distributions from the cases $n = 400$ and $n = 800$ are depicted in Figure 1. Clearly they are very close when $c_1$ is varying from 0 to $10^5$ for each case of $(n, h)$. As expected, they should depend on the bandwidth $h$. This suggests that the asymptotic null distribution of $l(\Theta_{0s}|G)$ is not very sensitive to the heteroscedasticity of the stochastic error. To check whether the scaled SELR statistics follow asymptotically the $\chi^2$-distribution, we equate the mean and variance of the scaled SELR, $r_0 l(\Theta_{0s}|G)$, to the corresponding mean and variance of a chi-squared random variable, say $\chi^2_{d_0}$, with degrees of freedom $d_0$. This results in $r_0 = 2\mu/\sigma^2$ and $d_0 = 2\mu^2/\sigma^2$ with $\mu$ and $\sigma^2$ the simulated mean and variance of $l(\Theta_{0s}|G)$. We calculated further the empirical distribution of the scaled SELR and compared it with the $\chi^2_{d_0}$-distribution for each combination of $(n, h)$. Since the empirical distributions do not depend sensitively on the conditional variance function, only one of them was used for comparison. As an example, Figure 2 depicts the two distributions for the case that

$$(n, h) = (800, 1.5 \times 800^{-2/9}) \quad \text{and} \quad c_1 = 1.$$

They are indeed very close. This demonstrates empirically the accuracy of the approximation of the null distribution of the proposed SELR statistic by using the $\chi^2$-distribution. We also conducted a similar simulation study for testing homogeneity:

$$H_{0p} : a_1(\cdot) = \theta.$$

It again shows that the Wilks phenomenon continues to hold for some composite null hypothesis testing problem. The details are not reported here.

To conclude this section, the power functions of the proposed test of $H_0$ are estimated and compared to the commonly used $F$-type test statistic,

$$F_{0s} = (\text{RSS0} - \text{RSS1})/\text{RSS1}$$

[see, e.g., Fan, Zhang and Zhang (2001), page 155 for the definition], based on 100 simulations for the sample sizes $n = 200, 800$ under two sequences of alternatives indexed by $r$. One is

(4.1) $$H_1 : a_1(u) = r(u - 0.5), \qquad r = 0.1, 0.2, \ldots,$$

TABLE 1
*Summary of simulation results. $\mu$ and $\sigma$ are the simulated mean and standard deviation of the SELR statistic based on 100 repetitions. $n$ is the sample size and $h$ is the bandwidth*

| $n$ | $h$ | $\mu$ | $\sigma$ | $\mu$ | $\sigma$ | $\mu$ | $\sigma$ |
|---|---|---|---|---|---|---|---|
| | | \multicolumn{6}{c}{Conditional variances} | | | | | |
| | | \multicolumn{2}{c}{1} | \multicolumn{2}{c}{$1+4u^2$} | \multicolumn{2}{c}{$1+10u^2$} |
| 100 | 0.35938 | 1.868 | 1.221 | 2.161 | 1.788 | 1.954 | 1.252 |
| 100 | 0.53907 | 1.754 | 1.480 | 1.767 | 1.642 | 1.646 | 1.219 |
| | | \multicolumn{6}{c}{Conditional variances} | | | | | |
| | | \multicolumn{2}{c}{1} | \multicolumn{2}{c}{$1+10u^2$} | \multicolumn{2}{c}{$1+100u^2$} |
| 200 | 0.15404 | 4.371 | 2.495 | 4.393 | 2.468 | 3.973 | 2.334 |
| 200 | 0.30808 | 2.463 | 1.527 | 2.263 | 1.421 | 2.215 | 1.413 |
| 200 | 0.46212 | 1.655 | 1.124 | 1.698 | 1.130 | 1.329 | 0.840 |
| 200 | 0.61616 | 1.376 | 1.081 | 1.519 | 1.242 | 1.395 | 0.976 |
| 400 | 0.13205 | 5.019 | 2.019 | 4.487 | 1.977 | 4.459 | 1.968 |
| 400 | 0.26410 | 3.081 | 1.720 | 2.681 | 1.361 | 2.965 | 1.433 |
| 400 | 0.39615 | 2.007 | 1.271 | 1.961 | 1.246 | 2.192 | 1.307 |
| 400 | 0.52820 | 1.867 | 1.492 | 1.622 | 1.126 | 1.743 | 1.501 |
| | | \multicolumn{6}{c}{Conditional variances} | | | | | |
| | | \multicolumn{2}{c}{1} | \multicolumn{2}{c}{$1+u^2$} | \multicolumn{2}{c}{$1+10^5 u^2$} |
| 800 | 0.12452 | 5.080 | 1.774 | 4.950 | 1.611 | 4.807 | 1.557 |
| 800 | 0.22640 | 3.092 | 1.354 | 3.191 | 1.457 | 3.093 | 1.455 |
| 800 | 0.33959 | 2.220 | 1.171 | 2.103 | 1.165 | 2.038 | 1.080 |
| 800 | 0.45279 | 1.785 | 1.078 | 1.627 | 1.069 | 1.699 | 1.069 |

and the other is

$$(4.2) \qquad H_1 : a_1(u) = r(2\sin^2(2\pi u) - 1), \qquad r = 0.1, 0.2, \ldots.$$

TABLE 2
*Empirical sizes of SELR and F-type tests. The probabilities are computed based on 100 simulations; $h = n^{-2/9}$ for $n = 200$ and $h = 1.5 n^{-2/9}$ for $n = 800$*

| | | \multicolumn{4}{c}{Conditional variances} | | \multicolumn{4}{c}{Conditional variances} |
| | | 1 | $1+u^2$ | $1+10u^2$ | $1+100u^2$ | | 1 | $1+u^2$ | $1+10u^2$ | $1+100u^2$ |
| $n$ | $c_r$ | \multicolumn{4}{c}{Sizes of SELR test} | $c_r$ | \multicolumn{4}{c}{Sizes of $F$-type test} |
|---|---|---|---|---|---|---|---|---|---|---|
| 200 | 5.20 | 0.05 | 0.05 | 0.07 | 0.07 | 0.0705 | 0.05 | 0.05 | 0.09 | 0.09 |
| 200 | 4.47 | 0.08 | 0.09 | 0.09 | 0.09 | 0.0579 | 0.09 | 0.11 | 0.13 | 0.15 |
| 200 | 3.16 | 0.22 | 0.25 | 0.25 | 0.24 | 0.0375 | 0.25 | 0.28 | 0.34 | 0.35 |
| 800 | 5.11 | 0.02 | 0.02 | 0.02 | 0.01 | 0.0134 | 0.02 | 0.05 | 0.09 | 0.09 |
| 800 | 4.59 | 0.04 | 0.03 | 0.03 | 0.02 | 0.0132 | 0.03 | 0.05 | 0.09 | 0.09 |
| 800 | 3.65 | 0.09 | 0.09 | 0.09 | 0.08 | 0.0109 | 0.09 | 0.10 | 0.12 | 0.17 |
| 800 | 2.81 | 0.20 | 0.21 | 0.21 | 0.19 | 0.00776 | 0.20 | 0.22 | 0.27 | 0.29 |

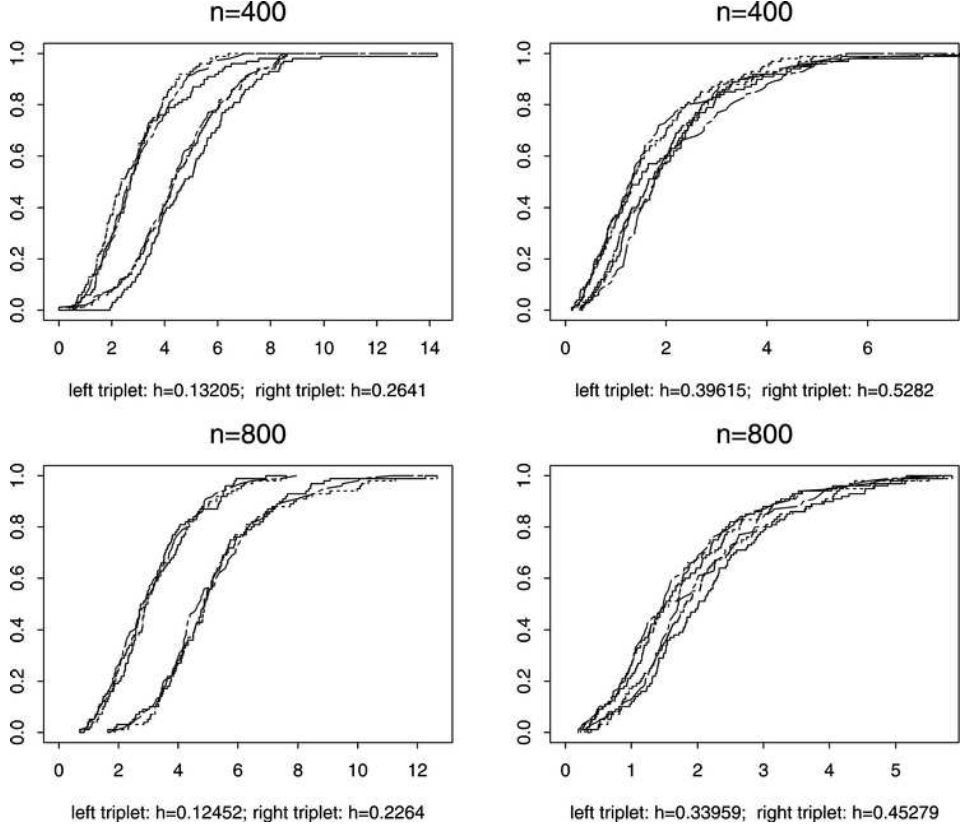

FIG. 1. *Comparisons of the empirical null distributions of the SELR statistics based on* 100 *simulations for different conditional variance functions. The solid curve, dotted curve, and dashed curve correspond to the conditional variances:* 1, $1 + 10u^2$ *and* $1 + 10^2 u^2$, *respectively, when* $n = 400$*; and to the conditional variances:* 1, $1 + u^2$ *and* $1 + 10^5 u^2$, *respectively, when* $n = 800$.

Here we take $1 + c_1 u^2$ as the conditional variance of the stochastic error given $U = u$ with $c_1 = 0, 1, 10, 10^2$. Note that the powers of the SELR and $F$-type test statistics have the same optimal rate $n^{-2/9}$. Thus, in this study for simplicity we select the bandwidth by comparing several empirically specified bandwidths. We find that the combinations of $h = n^{-2/9}$ and $n = 200$ and $h = 1.5 \times n^{-2/9}$ and $n = 800$ give relatively reasonable power functions for the two alternative sequences (4.1) and (4.2). For critical values given in Table 2, the sizes of the SELR and $F$-test are reported. It is evident that the sizes of the SELR test are adaptive automatically to the conditional variance function, while those of the $F$-type of test are not. This is consistent with our theoretical results and reflects one advantage of the SELR test.

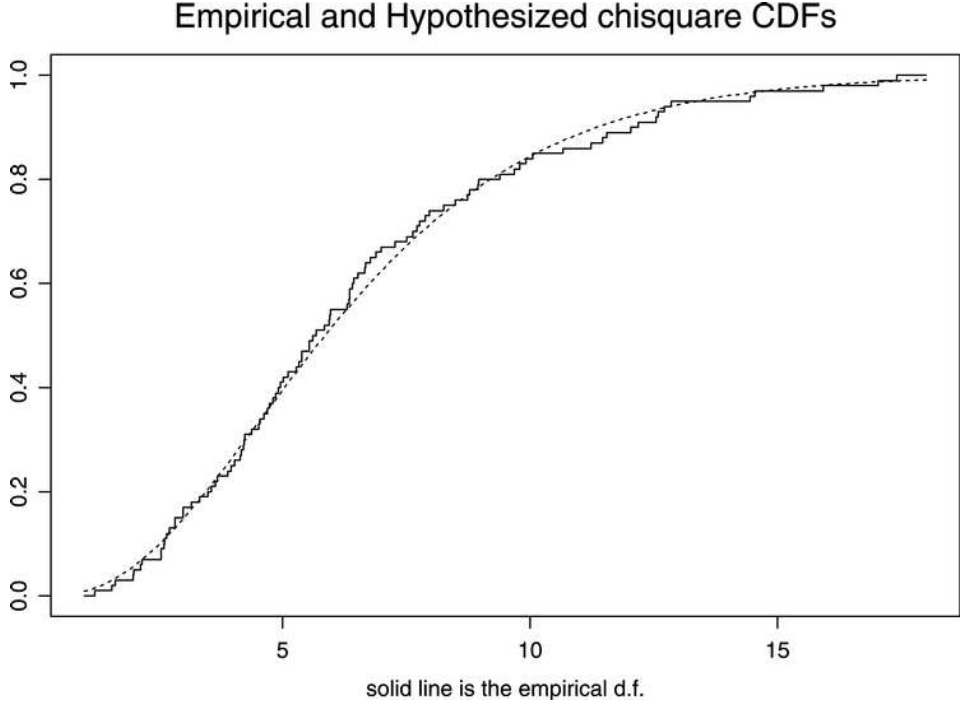

Fig. 2. *The solid stairstep curve is the empirical null distribution of the scaled sieve empirical likelihood ratio statistic for $n = 800$, $h = 1.5 \times 800^{-2/9}$, $c_1 = 1$ based on 100 repetitions, and the dashed curve is the chi-squared distribution with 6.51 degrees of freedom.*

Figures 3–5 present power functions at the significance levels shown in Table 2. We have conducted simulations on much more different settings and these are not reported to save space. As expected, the power deteriorates as the level of noise $c_1$ increases for both the SELR and $F$-type tests. Figures 3 and 5 indicate that the SELR test may significantly out-perform the $F$-type test in terms of power under the alternative

$$H_1 : a_1(u) = r(u - 0.5)$$

when there is heteroscedasticity. Similarly, Figure 4 implies that when the level of heteroscedasticity is low, the $F$-type test can have better power than the SELR test, and can perform much worse than the SELR test when the level of noise (heteroscedasticity) is high. This phenomenon can be explained by using Theorem 5. For the simple model $Y = a_1(U) + \varepsilon$, Theorem 5 gives

$$(4.3) \quad 2l(H_{0s}|G) = \frac{p}{h} K^*(0) + nE\{a_1(U)\sigma^{-2}(U)\} \\ + o(1) + (1 + o_p(h^{1/2}))\{T_n^* + 2W_{1n}^*\}.$$

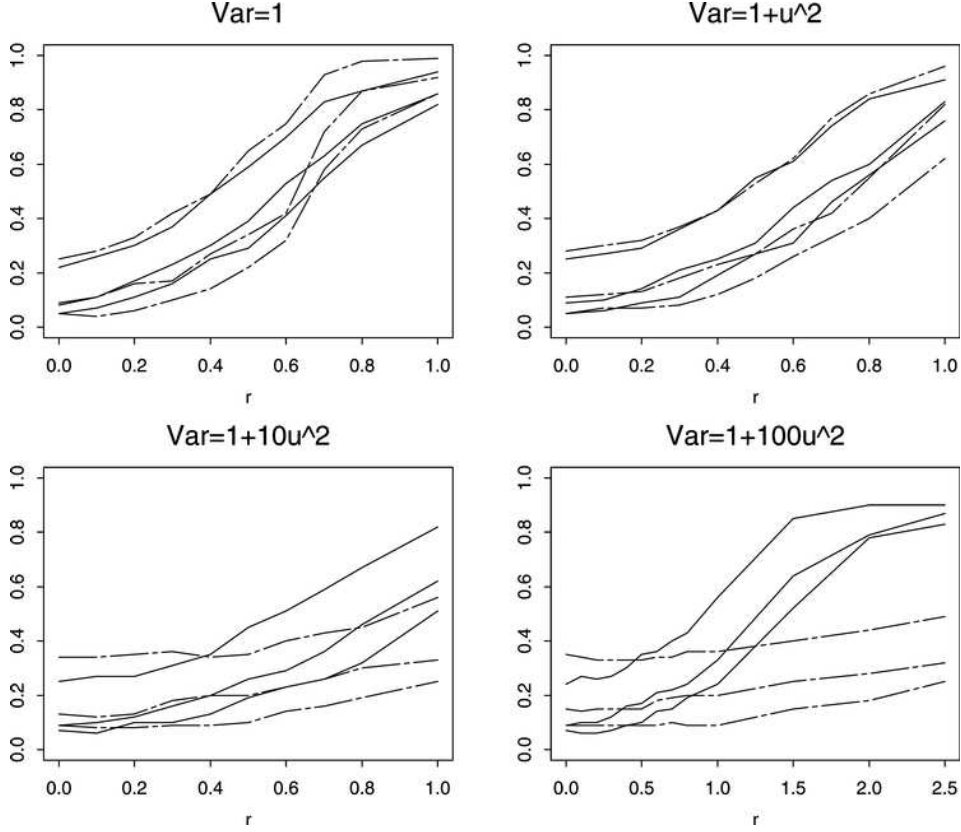

Fig. 3. *Comparisons of the power functions of the SELR (solid curves) and F-type tests (dashed curves) of $H_0 : a_1(u) = 0$ for $n = 200$ and the bandwidth $h = 200^{-2/9}$, evaluated at the alternatives* (4.1) *for different conditional variance functions.*

Note that if the function $\sigma(U)$ is known, we can make the transformation

$$Y' = Y/\sigma(U) + \varepsilon/\sigma(U)$$

and obtain the same asymptotic expansion as in (4.3) for the SELR based on the above transformed model. This means that (4.3) is adaptive to $\sigma^{-2}(U)$ in the sense that we can test $H_{0s} : a_1(\cdot) = 0$ asymptotically equally well whether or not we know the conditional variance of $\varepsilon$. In contrast, the $F$-type test does not have this property.

## 5. Technical conditions and proofs.

5.1. *Technical conditions.* Define

$$A_n(u_0, \beta) = \frac{1}{n} \sum_{i=1}^{n} K_h(U_i - u_0) \mathbf{G}_{ih}(u_0, \beta),$$

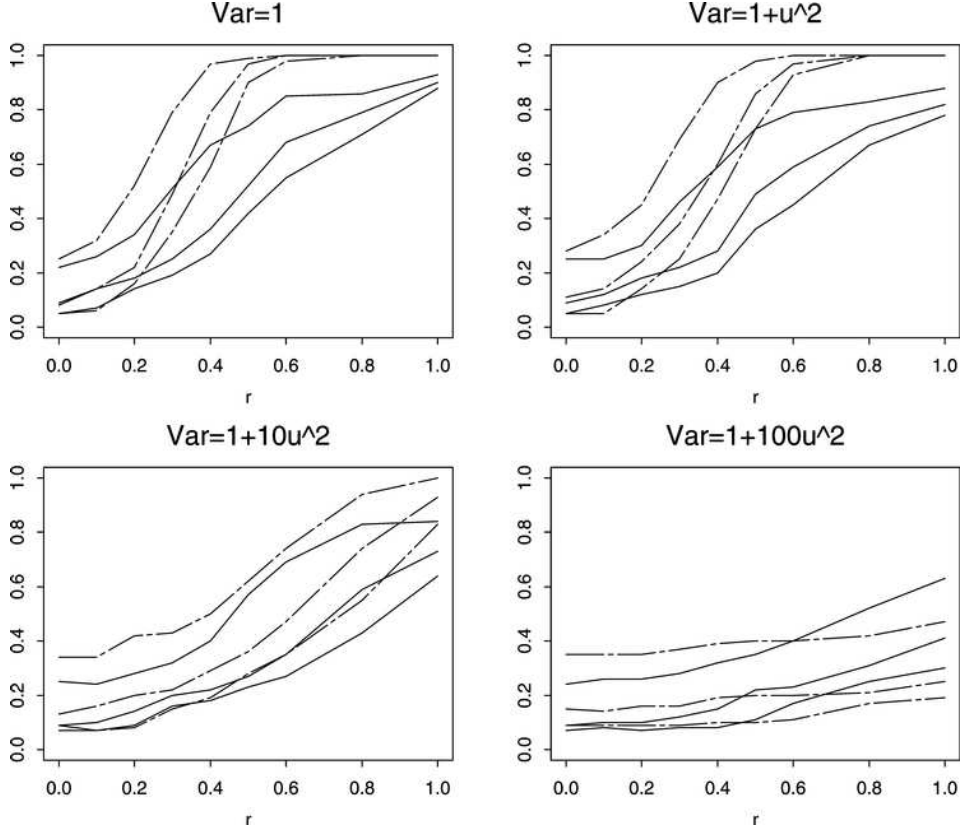

Fig. 4. *Comparisons of the power functions of the SELR (solid curves) and F-type tests (dashed curves) of $H_0: a_1(u) = 0$ for $n = 200$ and the bandwidth $h = 200^{-2/9}$, evaluated at the alternatives* (4.2) *for different conditional variance functions.*

$$Z_n(u_0, \beta) = \max_{1 \leq j \leq n} \|\mathbf{G}_{jh}(u_0, \beta)\|,$$

$$V_n(u_0, \beta) = \frac{1}{n} \sum_{i=1}^{n} K_h(U_i - u_0) \mathbf{G}_{ih}(u_0, \beta) \mathbf{G}_{ih}^{\tau}(u_0, \beta),$$

$$V_n(u_0, \alpha, \beta) = \frac{1}{n} \sum_{i=1}^{n} K_h(U_i - u_0) \frac{\mathbf{G}_{ih}(u_0, \beta) \mathbf{G}_{ih}^{\tau}(u_0, \beta)}{1 + \alpha^{\tau} \mathbf{G}_{ih}(u_0, \beta)},$$

$$B_n(u_0, \alpha, \beta) = \frac{1}{n} \sum_{i=1}^{n} \frac{K_h(U_i - u_0)}{1 + \alpha^{\tau} \mathbf{G}_{ih}(u_0, \beta)} \frac{\partial \mathbf{G}_{ih}(u_0, \beta)}{\partial \beta^{\tau}},$$

$$C_n(u_0, \alpha, \beta) = \frac{1}{n} \sum_{i=1}^{n} \frac{K_h(U_i - u_0)}{(1 + \alpha^{\tau} \mathbf{G}_{ih}(u_0, \beta))^2} \frac{\partial \mathbf{G}_{ih}(u_0, \beta)}{\partial \beta^{\tau}} \alpha \mathbf{G}_{ih}^{\tau}(u_0, \beta),$$

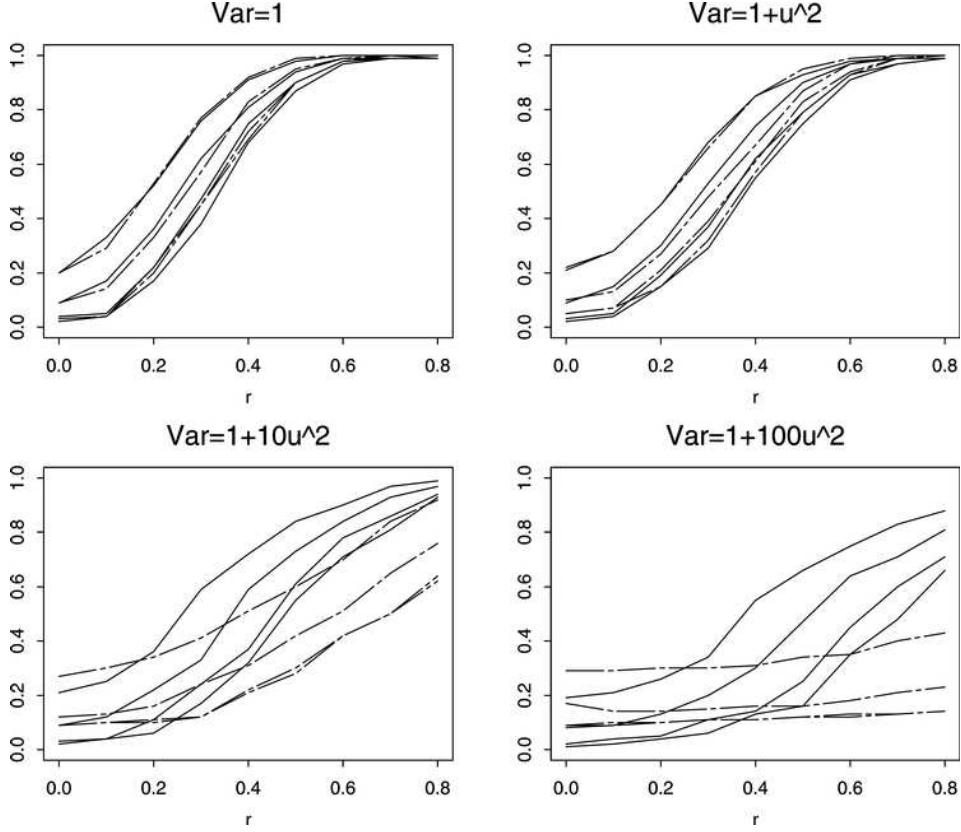

FIG. 5. *Comparisons of the power functions of the SELR (solid curves) and F-type tests (dashed curves) of $H_0: a_1(u) = 0$ for $n = 800$ and the bandwidth $h = 1.5 \times 800^{-2/9}$, evaluated at the alternatives* (4.1) *for different conditional variance functions.*

$$D_n(u_0, \alpha, \beta) = \frac{1}{n} \sum_{i=1}^{n} \frac{K_h(U_i - u_0)}{1 + \alpha^\tau \mathbf{G}_{ih}(u_0, \beta)} \frac{\partial^2 \mathbf{G}_{ih}(u_0, \beta)}{\partial \beta \partial \beta^\tau},$$

$$E_n(u_0, \alpha, \beta) = \frac{1}{n} \sum_{i=1}^{n} \frac{K_h(U_i - u_0)}{(1 + \alpha^\tau \mathbf{G}_{ih}(u_0, \beta))^2} \frac{\partial \mathbf{G}_{ih}(u_0, \beta)}{\partial \beta^\tau} \alpha \beta^\tau \frac{\partial \mathbf{G}_{ih}(u_0, \beta)^\tau}{\partial \beta}.$$

Here and hereafter the norm of a matrix $W = (w_{ij})$ is defined by $\|W\| = \sqrt{\sum_{i,j}^{n} w_{ij}^2}$. Let $r_0$ denote an arbitrary positive constant. Let $\Theta_0$ be a compact subset of $R^{2p}$ such that $\beta_0$ is an inner point of $\Theta_0$. Define

$$\mathcal{F}_0 = \{K((\cdot - u_0)/h) I\{\mathbf{G}_h(u_0, \beta)^\tau \psi > \delta\} : u_0 \in \Omega, \|\beta - \beta_0\| \leq r_0,$$
$$\|\psi\| = 1, 0 \leq \delta \leq 1\},$$

where $I\{\cdot\}$ is the indicator function,

$$\mathcal{F}_1 = \{K((\cdot - u_0)/h)\mathbf{G}_h(u_0, \beta) : u_0 \in \Omega, \|\beta - \beta_0\| \leq r_0\},$$
$$\mathcal{F}_2 = \{K((\cdot - u_0)/h)\mathbf{G}_h(u_0, \beta)\mathbf{G}_h^\tau(u_0, \beta)] : u_0 \in \Omega, \|\beta - \beta_0\| \leq r_0\},$$
$$\mathcal{F}_3 = \left\{K((\cdot - u_0)/h)\frac{\partial \mathbf{G}_h(u_0, \beta)}{\partial \beta^\tau} : u_0 \in \Omega, \beta \in \Theta_0\right\}.$$

Let $P_n$ denote the empirical distribution of $\{(U_i, X_i, Y_i)\}$, and $N(\delta, L_1(P_n), \mathcal{F}_j)$, $j = 0, 1, 2, 3$, the covering numbers [see, e.g., Pollard (1984), page 25 for the definition]. We impose the following technical conditions:

- (K0) $K$ has support $[-1, 1]$ and $\max_t K(t) < \infty$.
- (U0) The density of $U$ is Lipschitz continuous and bounded away from zero.
- (A1) $E[G(\varepsilon)|U] = 0$ and $\varepsilon$ is independent of $X$ given $U$.
- (A2) There exist a constant $\xi \geq 4$ and a function $F(Y, X)$ satisfying

$$\sup_{\substack{|t| \leq 1 \\ \|\beta - \beta_0\| \leq \delta_0}} \|G(Y - \beta^\tau Z(X, t))\|\|Z(X, t)\| \leq F(Y, X),$$

$$\sup_u E[F(Y, X)^\xi | U = u] < \infty.$$

- (A3) For $1 \leq k \leq k_0$,

$$\sup_{\substack{\|\beta - \beta_0\| \leq r_0 \\ u_0 \in \Omega, |t| \leq 1}} E[G_k^2(Y - \beta^\tau Z(X, t))\|Z(X, t)\|^2 | U = u_0 + th] = O(1).$$

- (A4) There exist $c_0(P_n)$ and some positive constants $c_0$ and $w_0$ such that $Ec_0(P_n) \to c_0$ and

$$N(\delta, L_1(P_n), \mathcal{F}_0) \leq c_0(P_n)(h\delta)^{-w_0}.$$

There exist $c_1(P_n)$ and some positive constants $c_1$ and $w_1$ such that $Ec_1(P_n) \to c_1$ and

$$N(\delta, L_1(P_n), \mathcal{F}_1) \leq c_1(P_n)(h\delta)^{-w_1}.$$

- (A5) Uniformly for $u_0 \in \Omega$, $\|t\| \leq 1$, $\|\beta - \beta_0\| \to 0$ and $h \to 0$,

$$E\left\{G\left(Y - \beta^\tau Z\left(X, \frac{U - u_0}{h}\right)\right)\Big| U = u_0 + th\right\} = O(h^2) + O(\|\beta - \beta_0\|).$$

- (A6) There exist $c_2(P_n)$ and a positive constant $c_2$ such that $Ec_2(P_n) \to c_2$ and

$$N(\delta, L_1(P_n), \mathcal{F}_2) \leq c_2(P_n)(h\delta)^{-w_2}.$$

(A7) $\sup_{\|\beta-\beta_0\|\leq r_0, u_0\in\Omega, |t|\leq 1} E[G_k^4(Y - \beta^\tau Z(X,t))\|Z(X,t)\|^4 | U = u_0 + th] = O(1)$.

(A8) Uniformly for $\|\beta - \beta_0\| \to 0$ and $h \to 0$,

$$E\left\{G\left(Y - \beta^\tau Z\left(X, \frac{U-u_0}{h}\right)\right)G^\tau\left(Y - \beta^\tau Z\left(X, \frac{U-u_0}{h}\right)\right)\Big|U\right\}$$
$$= V(u_0) + O(h^2) + O(\|\beta - \beta_0\|).$$

(A9) $V(u_0)$ and $\Gamma(u_0)$ defined in Section 3.1 are Lipschitz continuous in $u_0 \in \Omega$. Their minimum eigenvalues are uniformly positive in $u_0 \in \Omega$.

(A10) For any $\rho > 0$, there exists a constant $c(\rho) > 0$ such that when $h$ is small enough,

$$\inf_{\substack{\beta \in \Theta_0 \\ \|\beta - \beta_0\| \geq \rho}} \|EK_h(U - u_0)\mathbf{G}_h(u_0, \beta)\| > c(\rho).$$

For a positive sequence $\rho_{n1} \to 0$ and a small enough constant $\rho_2$, as $n \to \infty$,

$$\inf_{\rho_{n1} \leq \|\beta - \beta_0\| \leq \rho_2} \|EK_h(U - u_0)\mathbf{G}_h(u_0, \beta)\| \geq \rho_{n1} + O(h^2).$$

(B1) There exist a constant $\nu \geq 2$ and a function $F_4(Y, X)$ such that

$$\sup_u E[F_4^\nu(Y, X)|U = u] < \infty,$$

$$\sup_{u_0, \beta}\left\|\frac{\partial \mathbf{G}_h(u_0, \beta)}{\partial \beta^\tau}\right\| I(|U - u_0| \leq h) \leq F_4(Y, X).$$

(B2) For a constant $c$,

$$N(\delta, L_1(P_n), \mathcal{F}_3) \leq c(h\delta)^{-w_3}.$$

(B3) Uniformly for $u_0 \in \Omega$ and $\|\beta - \beta_0\| \leq r_n = o(h^{1/2})$,

$$E\left\{K_h(U - u_0)\frac{\partial \mathbf{G}_h(u_0, \beta)}{\partial \beta^\tau}\right\} = D(u_0) \otimes (S \otimes \Gamma(u_0)) + o(h^{1/2}).$$

(B4) $\sup_{\|\beta-\beta_0\|\leq r_0, u_0\in\Omega, |t|\leq 1} E[\|\partial \mathbf{G}_h(u_0, \beta)/\partial \beta^\tau\|^2 | U = u_0 + th] < \infty$.

(B5) There exists a function $F_5(y, x)$ such that

$$\sup_u E[F_5^2(Y, X)|Y = u] < \infty,$$

$$\sup_{\substack{u_0 \in \Omega \\ \|\beta - \beta_0\| \leq r_0}} \left\|\frac{\partial^2 \mathbf{G}_h(u_0, \beta)}{\partial \beta \partial \beta^\tau}\right\| I(|U - u_0| \leq h) \leq F_5(Y, X).$$

(B6) There exists a function $F_6$ such that $\sup_u E[F_6(\varepsilon, X) \|X\|^2 | U = u] < \infty$, and that for $|U - u_0| \leq h$ and

$$\varepsilon^* = \varepsilon + \frac{h^2}{2} A''^\tau(u_0 + s(U - u_0)) \frac{X(U - u_0)^2}{h^2}$$
$$+ (\beta - \beta_0)^\tau Z\left(X, \frac{U - u_0}{h}\right)$$

we have

$$\left\|\frac{\partial G(\varepsilon^*)}{\partial \varepsilon}\right\| \leq F_6(\varepsilon, X)$$

uniformly for $|s| \leq 1$, $\|\beta - \beta_0\| \leq r_0$ and $u_0 \in \Omega$.

We would like to make some comments on the conditions above. Suppressing dependence on $X$, we denote $Z(t) = Z(X, t)$. Suppose for some $r_0 > 0$ there exist integrable functions $F_j(Y, X), j = 1, 2, 3$, such that

$$\sup_{\|\beta - \beta_0\| \leq r_0, t} K'(t) \|G(Y - \beta^\tau Z(t))\| \|Z(t)\| \leq F_1(Y, X),$$

$$\sup_{\|\beta - \beta_0\| \leq r_0, t} K(t) \left\|\frac{\partial G(Y - \beta^\tau Z(t))}{\partial \varepsilon}\right\| \|Z(t)\|(\|Z'(t)\| + \|Z(t)\|) \leq F_2(Y, X),$$

$$\sup_{\|\beta - \beta_0\| \leq r_0, t} K(t) \|G(Y - \beta^\tau Z(t))\| \|Z'(t)\| \leq F_3(Y, X).$$

Then for some positive constant $c$,

$$\left\|K\left(\frac{u - u_1}{h}\right) \mathbf{G}_h(u_1, \beta_1) - K\left(\frac{u - u_2}{h}\right) \mathbf{G}_h(u_2, \beta_2)\right\|$$
$$\leq c\{F_1(Y, X) + F_2(Y, X) + F_3(Y, X)\}\left\{\frac{|u_1 - u_2|}{h} + \|\beta_1 - \beta_2\|\right\}.$$

Thus the second part of condition (A4) holds if $EF_j(Y, X) < \infty, j = 1, 2, 3$. Similar remarks can be made about conditions (A6) and (B2).

As pointed out in Section 2, $EK_h(U - u_0)\mathbf{G}_h(u_0, \beta_A) = 0, u_0 \in \Omega$, can be viewed as certain local estimating equations associated with the equations $E[G(Y - A(U)^\tau X) | U = u_0] = 0, u_0 \in \Omega$, as $A(u)$ is expanded around each $u_0$. In this sense, the first part of (A10) implies that when $\beta_A$ (coefficients of the approximation of $A$) is away from the true value $\beta_0$ (coefficients of the approximation of $A_0$), $\|EK_h(U - u_0)\mathbf{G}_h(u_0, \beta)\|$ is away from 0. This is a little stronger than the requirement that $E[G(Y - A^\tau(U)X) | U] = 0$ if and only if $A$ is equal to the true value. The second part of (A10) is a local condition which says locally $\|EK_h(U - u_0)\mathbf{G}_h(u_0, \beta)\|$ is bounded below by the norm of the linear function of $\beta$ near the true value $\beta_0$. For instance, assume the first component of $G$ is $Y - A^\tau(U)X$ and assume that $E[XX^\tau|U=u]$ is positive definite uniformly in $u$. Then we have

$$\|EK_h(U-u_0)\mathbf{G}_h(u_0,\beta_A)\|$$
$$\geq \left\|EK_h(U-u_0)\left[Y - \beta_A^\tau Z\left(X,\frac{U-u_0}{h}\right)\right] \otimes Z\left(X,\frac{U-u_0}{h}\right)\right\|$$
$$= O(h^2) + (\beta_0 - \beta_A)^\tau$$
$$\times \int K(t) E[Z(X,t)Z^\tau(X,t)|U=u_0+th] f(u_0+th)\,dt$$
$$\geq c\|\beta_0 - \beta_A\| + O(h^2),$$

provided $h$ is small enough.

5.2. *Proofs.* Note that Lemmas 1–8 are used in this section and their proofs can be found in the Appendix.

PROOF OF THEOREM 1. First of all, using Lemma 3, we obtain

$$\hat{\beta}(u_0) - \beta_0 = o_p(h^{1/2} \wedge n^{-1/\xi}), \qquad \hat{\alpha}(u_0) = o_p(h^{1/2} \wedge n^{-1/\xi}).$$

Furthermore, by the definition of $\hat{\alpha}\ (= \hat{\alpha}(u_0))$ and $\hat{\beta}\ (= \hat{\beta}(u_0))$, we have

$$0 = \frac{1}{n}\sum_{i=1}^n K_h(U_i - u_0) \frac{\mathbf{G}_{ih}(u_0,\hat{\beta})}{1 + \hat{\alpha}^\tau \mathbf{G}_{ih}(u_0,\hat{\beta})},$$
$$0 = \frac{1}{n}\sum_{i=1}^n K_h(U_i - u_0) \frac{\hat{\alpha}^\tau \partial \mathbf{G}_{ih}(u_0,\hat{\beta})/\partial \beta^\tau}{1 + \hat{\alpha}^\tau \mathbf{G}_{ih}(u_0,\hat{\beta})}.$$

Then invoking the Taylor expansion we have

$$0 = A_n(u_0, \beta_0) - V_n(u_0, \alpha_{n1}, \beta_{n1})\hat{\alpha}$$
$$\qquad + \{B_n(u_0, \alpha_{n1}, \beta_{n1}) - C_n(u_0, \alpha_{n1}, \beta_{n1})\}(\hat{\beta} - \beta_0),$$
$$0 = \{B_n(u_0, \alpha_{n2}, \beta_{n2}) - C_n(u_0, \alpha_{n2}, \beta_{n2})\}\hat{\alpha}$$
$$\qquad + \{D_n(u_0, \alpha_{n2}, \beta_{n2}) - E_n(u_0, \alpha_{n2}, \beta_{n2})\}(\hat{\beta} - \beta_0),$$

where $\alpha_{nj}$, $j=1,2$, are between $\hat{\alpha}$ and $0$ and $\beta_{nj}$, $j=1,2$, are between $\hat{\beta}$ and $\beta_0$. By using Lemmas 4–8, the above equations become

$$-A_n(u_0, \beta_0) = -(1 + o_p(h^{1/2}))V(u_0) \otimes (S \otimes \Gamma(u_0))\hat{\alpha}$$
$$\qquad + \{o_p(h^{1/2}) + D(u_0) \otimes (S \otimes \Gamma(u_0))\}(\hat{\beta} - \beta_0),$$
$$0 = \{o_p(h^{1/2}) + D(u_0) \otimes (S \otimes \Gamma(u_0))\}\hat{\alpha} + o_p(h^{1/2})(\hat{\beta} - \beta_0).$$

It follows that

$$(\hat{\beta} - \beta_0) = -[(D(u_0)^\tau V^{-1}(u_0) D(u_0))^{-1} D(u_0) V^{-1}(u_0)$$
$$\otimes (S^{-1} \otimes \Gamma(u_0)^{-1}) + o_p(h^{1/2})] A_n(u_0, \beta_0),$$
$$\hat{\alpha} = [V^{-1}(u_0) - V^{-1}(u_0)(D^\tau(u_0) V^{-1}(u_0) D(u_0))^{-1}$$
$$\times D(u_0) D^\tau(u_0) V^{-1}(u_0) + o_p(h^{1/2})] A_n(u_0, \beta_0).$$

Observe that for $U_i^* = u_0 + s(U_i - u_0)$, $0 \leq s \leq 1$, and for

$$\varepsilon_i^* = Y_i - A^\tau(U_i) X_i + \frac{1}{2} A''(U_i^*) X_i (U_i - u_0)^2$$
$$+ (\beta - \beta_0)^\tau Z\left(X_i, \frac{U_i - u_0}{h}\right)$$

we have

$$A_n(u_0, \beta) = \frac{1}{n} \sum_{i=1}^n K_h(U_i - u_0) G(\varepsilon_i^*) \otimes Z\left(X_i, \frac{U_i - u_0}{h}\right)$$
$$= \frac{1}{n} \sum_{i=1}^n K_h(U_i - u_0) G(\varepsilon_i) \otimes Z\left(X_i, \frac{U_i - u_0}{h}\right)$$
$$+ \frac{h^2}{2} O_p(1) + O_p(\|\beta - \beta_0\|),$$

where the last equality follows from the condition (B6) (or $A$ is linear). Now the proof can be completed by some simple calculations. $\square$

PROOF OF THEOREM 2. Note that under the conditions of Theorem 1 we have $h \to 0$ and $nh^{3/2} \to \infty$. Recall that given $U$, $\varepsilon$ and $X$ are independent by condition (A1). By the Taylor expansion and Lemma 4 there are matrices $V_n^*(U_j)$ such that as $n \to \infty$, uniformly in $U_j$,

$$V_n^*(U_j) = V(U_j) \otimes (S \otimes \Gamma(U_j))(1 + o_p(h^{1/2})),$$
$$\hat{\alpha}(U_j) = V_n^*(U_j)^{-1} \frac{1}{n} \sum_{i=1}^n K_h(U_i - U_j) \mathbf{G}_{ih}(U_j, \hat{\beta}).$$

The last two equalities lead to

$$l(G) = \sum_{j=1}^n \hat{\alpha}(U_j)^\tau \sum_{i=1}^n \frac{K_h(U_i - U_j)}{\sum_{m=1}^n K_h(U_m - U_j)} \mathbf{G}_{ih}(U_j, \hat{\beta})$$
$$- \frac{1}{2} \sum_{j=1}^n \hat{\alpha}(U_j)^\tau V_n(U_j, s^* \hat{\alpha}(U_j), \hat{\beta}) \hat{\alpha}(U_j)$$

$$(5.1) \quad = \sum_{j=1}^{n} \hat{\alpha}(U_j)^\tau \left\{ \frac{1}{\sum_{m=1}^{n} K_h(U_m - U_j)} V_n^*(U_j) \right.$$

$$\left. - \frac{1}{2} V_n(U_j, s^* \hat{\alpha}(U_j, \hat{\beta})) \right\} \hat{\alpha}(U_j)$$

$$= \frac{1}{2}(1 + o_p(h^{1/2})) \sum_{j=1}^{n} f^{-1}(U_j) \hat{\alpha}(U_j)^\tau [V(U_j) \otimes (S \otimes \Gamma(u_0))] \hat{\alpha}(U_j),$$

where $0 \leq s^* \leq 1$, and $V_n(u, \alpha, \beta)$ is defined in Section 5.1. Note that we draw out the factor $1 + o_p(h^{1/2})$ from the inside of the summation in (5.1) because the $o_p(h^{1/2})$ is uniform with respect to $U_j$, $1 \leq j \leq n$, and $\hat{\alpha}(U_j)^\tau[V(U_j) \otimes (S \otimes \Gamma(U_j))]\hat{\alpha}(U_j)/\hat{\alpha}(U_j)^\tau \hat{\alpha}(U_j)$, $1 \leq j \leq n$, are bounded away from 0 and $\infty$ [see condition (A9)]. It follows from the definition of $C(u)$ in Section 3.1 that $C(u)V(u)C(u) = C(u)$. Thus, combining (5.1) and Theorem 1, we obtain

$$l(G) = \left(\frac{1}{2} + o_p(h^{1/2})\right) \sum_{j=1}^{n} \frac{1}{n} \sum_{i=1}^{n} K_h(U_i - U_j)(C(U_j)G(\varepsilon_i))^\tau$$

$$\otimes \begin{pmatrix} \Gamma^{-1}(U_j)X_i \\ \mu_2^{-1}(U_i - U_j)\Gamma^{-1}(U_j)X_i/h \end{pmatrix}^\tau$$

$$\times \frac{1}{f(U_j)}[V(U_j) \otimes (S \otimes \Gamma(U_j))]$$

$$\times \frac{1}{n} \sum_{k=1}^{n} K_h(U_k - U_j)(C(U_j)G(\varepsilon_k))$$

$$(5.2) \qquad \otimes \begin{pmatrix} \Gamma^{-1}(U_j)X_k \\ \mu_2^{-1}(U_k - U_j)\Gamma^{-1}(U_j)X_k/h \end{pmatrix} + \zeta_n$$

$$= (1 + o_p(h^{1/2}))$$

$$\times \frac{1}{2n^2} \sum_{i=1}^{n} \sum_{k=1}^{n} \sum_{j=1}^{n} K_h(U_i - U_j) K_h(U_k - U_j) f^{-1}(U_j)$$

$$\times G^\tau(\varepsilon_i) C(U_j) G(\varepsilon_k)$$

$$\times \left(1 + \mu_2^{-1} \frac{(U_i - U_j)(U_k - U_j)}{h^2}\right)$$

$$\times X_i^\tau \Gamma^{-1}(U_j) X_k + \zeta_n,$$

where $\zeta_n = 0$ when $A_0$ is linear, and otherwise $\zeta_n = O_p(nh^4)$. The last term in (5.2) can be decomposed as follows:

$$(5.3) \quad (1 + o_p(h^{1/2}))L(G) = T_{n11} + T_{n121} + T_{n122} + T_{n21} + T_{n22} + \zeta_n$$

where

$$T_{n11} = \frac{1}{n^2} \sum_{i=1}^{n} \sum_{j=1}^{n} K_h(U_i - U_j)^2 f^{-1}(U_j)$$

$$\times \{G^\tau(\varepsilon_i)C(U_j)G(\varepsilon_i) - E[G^\tau(\varepsilon_i)C(U_j)G(\varepsilon_i)|(U_i, U_j)]\}$$

$$\times \left(1 + \mu_2^{-1}\frac{(U_i - U_j)^2}{h^2}\right) X_i^\tau \Gamma^{-1}(U_j) X_i,$$

$$T_{n121} = \frac{1}{n^2} \sum_{i=1}^n \sum_{j=1}^n K_h(U_i - U_j)^2 E[G^\tau(\varepsilon_i) C(U_j) G(\varepsilon_i) | (U_i, U_j)]$$

$$\times \left(1 + \mu_2^{-1}\frac{(U_i - U_j)^2}{h^2}\right)$$

$$\times \{X_i^\tau \Gamma^{-1}(U_j) X_i - E[X_i^\tau \Gamma^{-1}(U_j) X_i | (U_i, U_j)]\} f^{-1}(U_j),$$

$$T_{n122} = \frac{1}{n^2} \sum_{i=1}^n \sum_{j=1}^n K_h(U_i - U_j)^2 E[G^\tau(\varepsilon_i) C(U_j) G(\varepsilon_i) | (U_i, U_j)]$$

$$\times \left(1 + \mu_2^{-1}\frac{(U_i - U_j)^2}{h^2}\right)$$

$$\times E[X_i^\tau \Gamma^{-1}(U_j) X_i | (U_i, U_j)] f^{-1}(U_j),$$

$$T_{n21} = \frac{1}{n^2} \sum_{i \neq k} \sum_{j \notin \{i,k\}} K_h(U_i - U_j) K_h(U_k - U_j) G^\tau(\varepsilon_i) C(U_j) G(\varepsilon_k)$$

$$\times \left(1 + \frac{(U_i - U_j)(U_k - U_j)\mu_2^{-1}}{h}\right) X_i^\tau \Gamma^{-1}(U_j) X_k,$$

$$T_{n22} = \frac{K(0)}{n^2 h^2} \sum_{i \neq k}^n \bigg\{ K\left(\frac{U_k - U_i}{h}\right) G^\tau(\varepsilon_i) C(U_j) G(\varepsilon_k) X_i^\tau \Gamma^{-1}(U_i) X_k f^{-1}(U_i)$$

$$+ K\left(\frac{U_i - U_k}{h}\right) G^\tau(\varepsilon_i) C(U_k) G(\varepsilon_k) X_i^\tau \Gamma^{-1}(U_k) X_k f^{-1}(U_k) \bigg\}.$$

Observe that as $nh^{3/2} \to \infty$, $h \to 0$,

$$\begin{aligned}
T_{n122} &= \frac{K(0)^2}{(nh)^2} \sum_{i=1}^n \operatorname{tr}(C(U_i) V(U_i)) p f^{-2}(U_i) \\
&\quad + \frac{1}{n^2} \sum_{i \neq j} K_h(U_i - U_j)^2 \operatorname{tr}(C(U_j) V(U_i)) \\
&\quad \times \left(1 + \mu_2^{-1}\frac{(U_i - U_j)^2}{h^2}\right) \operatorname{tr}(\Gamma^{-1}(U_j) \Gamma(U_i))(f(U_i) f(U_j))^{-1} \\
&= \frac{K(0)^2}{nh^2} \bigg\{ E\bigg[\frac{\operatorname{tr}(C(U) V(U))}{f^2(U)}\bigg] + O_p(n^{-1/2}) \bigg\} + \Psi_n \\
&= o_p(h^{-1/2}) + \Psi_n
\end{aligned} \quad (5.4)$$

where

$$\Psi_n = \frac{1}{n^2} \sum_{i \neq j} K_h(U_i - U_j) \operatorname{tr}(C(U_j)V(U_i))$$

$$\times \left(1 + \frac{\mu_2^{-1}(U_i - U_j)^2}{h^2}\right) \operatorname{tr}(\Gamma^{-1}(U_j)\Gamma(U_i))(f(U_i)f(U_j))^{-1}$$

$$= \frac{(k_0 - 1)p|\Omega|}{h} \int K^2(t)(1 + \mu_2^{-1}t^2)\, dt + o_p(h^{-1/2}).$$

This is because

$$E\Psi_n = (1 + O(h))\frac{p}{h} \int K(t)^2(1 + \mu_2^{-1}t^2)\, dt\, E\{\operatorname{tr}(C(U)V(U))f^{-1}(U)\}$$

$$= \frac{p(k_0 - 1)}{h}(1 + O(h))|\Omega| \int K^2(t)(1 + \mu_2^{-1}t^2)\, dt,$$

$$\operatorname{Var}(\Psi_n) \leq O(n^{-1}h^{-2}) = o(h^{-1}).$$

By a similar argument, we have the following equalities:

$$T_{n121} = \frac{K(0)^2}{(nh)^2} \sum_{i=1}^n \operatorname{tr}(C(U_i)V(U_i))$$

$$\times (X_i^\tau \Gamma^{-1}(U_i)X_i - E[X_i^\tau \Gamma^{-1}(U_i)X_i|U_i])f^{-1}(U_i)$$

(5.5)
$$+ \frac{1}{n^2} \sum_{i \neq j} K_h(U_i - U_j)^2 \operatorname{tr}(C(U_j)V(U_j))\left(1 + \mu_2^{-1}\frac{(U_i - U_j)^2}{h^2}\right)$$

$$\times \{X_i^\tau \Gamma^{-1}(U_j)X_i - E[X_i^\tau \Gamma^{-1}(U_j)X_i|(U_i, U_j)]\}f^{-1}(U_j)$$

$$= \frac{K(0)^2}{nh^2}O_p(n^{-1/2}) + o_p(h^{-1/2}),$$

(5.6) $$T_{n22} = o_p(h^{-1/2}),$$

(5.7) $$T_{n21} = o_p(h^{-1/2}) + \frac{n-2}{n^2}\sum_{i \neq k} G^\tau(\varepsilon_i)\Phi_{ikh}G(\varepsilon_k),$$

where $\Phi_{ikh}$ is defined in (3.2) and the last equality follows from Hoeffding's decomposition for the variance of $U$-statistics. Now (5.3)–(5.7) imply (3.3). Equation (3.4) can be proved by a similar argument by showing that

$$l(A_0|G) = (1 + o_p(h^{1/2}))$$

$$\times \frac{1}{2n^2}\sum_{j=1}^n A_n^\tau(U_j, \beta_0)[V(U_j) \otimes (S \otimes \Gamma(U_j))]^{-1} A_n(U_j, \beta_0).$$

The proof is complete. □

PROOF OF THEOREM 3. Invoking the asymptotic representations in Theorem 2, we need only to prove the asymptotic normality of $T_n$. To this end, we first calculate the variance of $T_n$,

$$\begin{aligned}
\operatorname{Var}(T_n) &= \frac{(2+o(1))}{n(n-1)} E\{G^\tau(\varepsilon_1)\Phi_{12h}G(\varepsilon_2)\}^2 \\
&= \frac{(2+o(1))}{n(n-1)} \operatorname{tr}\{E[\Phi_{12h}G(\varepsilon_2)G^\tau(\varepsilon_2)\Phi_{123h}^\tau G(\varepsilon_1)G^\tau(\varepsilon_1)]\} \\
&= \frac{2(1+O(h))}{n(n-1)} \operatorname{tr}\{E[K_h^*(U_2-U_1)^2 C(U_1)G(\varepsilon_2)G^\tau(\varepsilon_2)C(U_1)G(\varepsilon_2) \\
&\qquad\qquad\qquad \times G(\varepsilon_2)^\tau X_1^\tau \Gamma^{-1}(U_1) X_2 X_2^\tau \Gamma^{-1}(U_1) X_1]\} \\
&= \frac{2(1+O(h))}{n(n-1)} \frac{p(k_0-1)}{h} |\Omega| \int K^*(t)^2\, dt.
\end{aligned}$$

(5.8)

Let $D_i = (\varepsilon_i, X_i, U_i)$, $1 \le i \le n$, and $\Pi_k$ be the $\sigma$-algebra generated by $D_1, \ldots, D_k$, $1 \le k \le n$. Set $\Phi_h(D_i, D_k) = G^\tau(\varepsilon_i)\Phi_{ikh}G(\varepsilon_k)$, $\eta_{n1}=0$ and

$$\eta_{nk} = E[T_n|\Pi_k] - E[T_n|\Pi_{k-1}].$$

Then

$$\eta_{nk} = \frac{2}{n(n-1)} \sum_{j=1}^{k-1} \Phi_h(D_j, D_k), \qquad 2 \le k \le n,$$

and $\{\eta_{nk}, \Pi_k\}$ is a sequence of martingale differences. By Theorem 4 of Shiryayev [(1996), page 543], it suffices to show

(5.9) $$\operatorname{Var}^{-1}(T_n) \sum_{k=2}^n E[\eta_{nk}^2|\Pi_{k-1}] \to 1 \qquad \text{in probability}$$

and

(5.10) $$\operatorname{Var}^{-2}(T_n) \sum_{k=1}^n E\eta_{nk}^4 \to 0.$$

In the following, $D = (\varepsilon, X, U)$ denotes a general random variable independent of $D_i$ and $D_k$. To prove (5.9) and (5.10), we need the following equalities for $i < j$:

$$\begin{aligned}
&E[\Phi_h(D_i, D_j)^2|D_i] \\
&\quad = \frac{1}{h} \int K^*(t)^2\, dt\, X_i^\tau \Gamma^{-1}(U_i) X_i G^\tau(\varepsilon_i) C(U_i) G(\varepsilon_i)(1+O(h)), \\
&E[\Phi_h(D_i, D)\Phi_h(D_j, D)|(D_i, D_j)] \\
&\quad = G(\varepsilon_i)^\tau E[K_h(U-U_i)K_h(U-U_k)C(U_i)V(U)C(U_j) \\
&\qquad\qquad \times \operatorname{tr}(\Gamma^{-1}(U_k)XX^\tau \Gamma^{-1}(U_j) X_i X_j^\tau)|(D_i, D_j)]G(\varepsilon_j),
\end{aligned}$$

$$E\Phi_h^2(D_i, D)\Phi_h^2(D_j, D)$$
$$= \frac{1}{h^2}(1 + O(h)) \iint K^*(t)^2 K^*(s)^2 \, dt \, ds$$
$$\times E[(X_j^\tau \Gamma^{-1}(U_j) X_j)^2 (G^\tau(\varepsilon_j) C(U_j) G(\varepsilon_j))^2],$$
$$E\Phi_h^4(D_i, D_k)$$
$$= O(1)(1 + O(h)) \frac{1}{h^3} \int K^*(t)^4 \, dt.$$

These are obvious by the assumption that $\varepsilon$ and $X$ are independent given $U$. Now with the above equalities, we can derive

$$\sum_{k=2}^n E[\eta_{nk}^2 | \Pi_{k-1}]$$

$$= \sum_{k=2}^n \frac{4}{n^2(n-1)^2} \left\{ \sum_{j=1}^{k-1} E[\Phi_h(D_j, D_k)^2 | D_j] \right.$$
$$\left. + \sum_{i \neq j}^{k-1} E[\Phi_h(D_i, D_k) \Phi_h(D_j, D_k) | (D_i, D_j)] \right\}$$

$$= \sum_{k=2}^n \frac{4}{n^2(n-1)^2} \sum_{i=1}^{k-1} \frac{1 + O(h)}{h} \int K^*(t)^2 \, dt$$

(5.11)
$$\times X_i^\tau \Gamma^{-1}(U_i) X_i G^\tau(\varepsilon_i) C(U_i) G(\varepsilon_i)$$

$$+ \sum_{k=2}^n \frac{8}{n^2(n-1)^2} \sum_{i<j}^{n-1} (n-j) E[\Phi_h(D_i, D) \Phi_h(D_j, D) | (D_i, D_j)]$$

$$= \frac{(1 + O(h)) 4 \int K^*(t)^2 \, dt}{n^2(n-1)^2 h}$$
$$\times \sum_{i=1}^{n-1} (n-i) X_i^\tau \Gamma^{-1}(U_i) X_i G^\tau(\varepsilon_i) C(U_i) G(\varepsilon_i) + \Upsilon_n$$

$$= (1 + o(1)) \frac{2 \int K^*(t)^2 \, dt}{n(n-1)}$$
$$\times E\{E[X_1^\tau \Gamma^{-1}(U_i) X_i | U_i] E[G^\tau(\varepsilon_i) C(U_i) G(\varepsilon_i) | U_i]\} + \Upsilon_n$$
$$= (1 + o(1)) \operatorname{Var}(T_n) + \Upsilon_n,$$

where

$$\Upsilon_n = \frac{8}{n^2(n-1)^2}$$

$$\times \sum_{i<j}^{n-1}(n-j)G(\varepsilon_i)^\tau$$
$$\times E[K_h^*(U-U_i)K_h^*(U-U_k)C(U_i)V(U)C(U_k)$$
$$\times \operatorname{tr}(\Gamma^{-1}(U_k)XX^\tau\Gamma^{-1}(U_i)X_iX_k^\tau)|(U_i,U_k,X_i,X_k)]G(\varepsilon_k).$$

Note that

$$E[\Upsilon_n]^2 = \frac{64}{n^4(n-1)^4}$$
$$\times \sum_{i<k}^{n-1}(n-k)^2$$
$$\times E\{G(\varepsilon_i)^\tau E[K_h^*(U-U_i)K_h^*(U-U_k)C(U_i)V(U)C(U_k)$$
$$\times \operatorname{tr}(\Gamma^{-1}(U_k)XX^\tau$$
$$\times \Gamma^{-1}(U_i)X_iX_k^\tau)|(U_i,X_i,U_k,X_k)]G(\varepsilon_k)\}^2$$
$$= O\left(\frac{(k_0-1)p}{(n-1)^4 h}\right)\int [K^*(t)*K^*(t)]^2\,dt$$
$$= O\left(\frac{1}{n^2}\right)\operatorname{Var}(T_n),$$

which implies $\Upsilon_n = o_p(\operatorname{Var}(T_n))$, and where $K^*(t)*K^*(t)$ is the convolution of $K^*(t)$ with itself. Substituting the above equality into (5.11), we get (5.9). Analogously, (5.10) follows from the following calculations:

$$\sum_{k=2}^n E\eta_{nk}^4 = \frac{O(1)}{n^4(n-1)^4}\sum_{k=2}^n\sum_{i\neq j}^{k-1}\left\{O\left(\frac{1}{h^2}\right)+O\left(\frac{k-1}{h^3}\right)\right\}$$
$$= O(\operatorname{Var}(T_n)^2)\frac{1}{(n-1)^2}\left\{O(n)+O\left(\frac{1}{h}\right)\right\}.$$

The proof is complete. □

PROOF OF THEOREM 4. The first part is similar to the proof of Theorem 3. The details are omitted. To show the second part, we recall that $\Gamma(u_0) = E[XX^\tau|U=u_0]f(u_0)$ and write

$$X_k = \begin{pmatrix} X_k^{(1)} \\ X_k^{(2)} \end{pmatrix}, \quad \Gamma = \begin{pmatrix} \Gamma_{11} & \Gamma_{12} \\ \Gamma_{21} & \Gamma_{22} \end{pmatrix} \quad \text{and} \quad \Gamma_{11,2} = \Gamma_{11} - \Gamma_{12}\Gamma_{22}^{-1}\Gamma_{21}$$

where $X_k^{(1)}$ is $p_1$-dimensional, $\Gamma_{11},\Gamma_{12},\Gamma_{21},\Gamma_{22}$ are $p_1\times p_1$, $p_1\times p_2$, $p_2\times p_1$ and $p_2\times p_2$ matrices and $p_2 = p-p_1$. Following the same steps in the proof of Theorem 3, we first extend Theorem 1 as follows:

$$\hat{\beta}_2(u_0) = \beta_2(u_0) + \frac{1}{n}\sum_{i=1}^{n} K_h(U_i - u_0) \begin{pmatrix} \Gamma_{22}^{-1}(u_0)X_i^{(2)} \\ \mu_2^{-1}\Gamma_{22}^{-1}(u_0)X_i^{(2)}(U_i - u_0)/h \end{pmatrix}$$
$$\times \eta_i(u_0)(1 + o_p(h^{1/2})) + O_p(h^2),$$

$$\hat{\alpha}^*(u_0) = \frac{1}{n}\sum_{i=1}^{n} K_h(U_i - u_0)$$
$$\times \left\{ (V^{-1}(u_0)G(\varepsilon_i)) \otimes \begin{pmatrix} \Gamma^{-1}(u_0)X_i \\ \mu_2^{-1}\Gamma^{-1}(u_0)X_i(U_i - u_0)/h \end{pmatrix} \right.$$
$$- V^{-1}D(D^\tau V^{-1}D)^{-1}D^\tau V^{-1}G(\varepsilon_i)$$
$$\left. \otimes \begin{pmatrix} 0 \\ \Gamma_{22}^{-1}(u_0)X_i^{(2)} \\ 0 \\ \mu_2^{-1}\Gamma_{22}^{-1}(u_0)X_i^{(2)}(U_i - u_0)/h \end{pmatrix} \right\}$$
$$\times (1 + o_p(h^{1/2})) + O_p(h^2).$$

Then by using the decomposition formula in Fan, Zhang and Zhang (2001) we have

$$X_i^\tau \Gamma(U_k)^{-1} X_k$$
$$= \{X_i^{(1)\tau} - X_i^{(2)\tau}\Gamma_{22}(U_k)\Gamma_{21}(U_k)\}\Gamma_{11.2}^{-1}(U_k)\{X_k^{(1)} - \Gamma_{12}(U_k)\Gamma_{22}^{-1}(U_k)X_k^{(2)}\}$$
$$+ X_i^{(2)\tau}\Gamma_{22}(U_k)^{-1}X_k^{(2)}.$$

The remaining part is very similar to the proof of Theorem 3. The details are omitted. □

PROOF OF THEOREM 5. The argument is similar to that in Fan, Zhang and Zhang (2001) but more tedious. For simplicity, we derive it heuristically. Write

$$l(H_{0s}|G) = (1 + o_p(h^{1/2}))$$
$$\times \frac{1}{2n^2}\sum_{i=1}^{n}\sum_{k=1}^{n}\sum_{j=1}^{n} K_h(U_i - U_j)K_h(U_k - U_j)$$
$$\times \frac{1}{f(U_j)}G^\tau(\varepsilon_i + A(U_i)^\tau X_i)(\varepsilon_k + A(U_k)^\tau X_k)$$

(5.12)

$$\times \left(1 + \mu_2^{-1}\frac{U_i - U_j}{h}\frac{U_k - U_j}{h}\right)$$
$$\times X_i^\tau \Gamma^{-1}(U_j) X_k - l_G$$
$$= (1 + o_p(h^{1/2}))(W_{n0} + W_{n1} + W_{n2} + W_{n3}) - l_G$$

with

$$W_{n0} = \frac{1}{2n^2} \sum_{i=1}^n \sum_{k=1}^n \sum_{j=1}^n K_h(U_i - U_j) K_h(U_k - U_j)$$
$$\times \left(1 + \mu_2^{-1}\frac{U_i - U_j}{h}\frac{U_k - U_j}{h}\right)\frac{1}{f(U_j)}$$
$$\times G(\varepsilon_i)^\tau V^{-1}(U_j) G(\varepsilon_k) X_i^\tau \Gamma^{-1}(U_j) X_k,$$

$$W_{n1} = \frac{1}{2n^2} \sum_{i=1}^n \sum_{k=1}^n \sum_{j=1}^n K_h(U_i - U_j) K_h(U_k - U_j)$$
$$\times \left(1 + \mu_2^{-1}\frac{U_i - U_j}{h}\frac{U_k - U_j}{h}\right)\frac{1}{f(U_j)}$$
$$\times G(\varepsilon_i)^\tau V^{-1}(U_j) \frac{\partial G(\varepsilon_k^*)}{\partial \varepsilon} X_i^\tau \Gamma^{-1}(U_j) X_k X_k^\tau A(U_k),$$

$$W_{n2} = \frac{1}{2n^2} \sum_{i=1}^n \sum_{k=1}^n \sum_{j=1}^n K_h(U_i - U_j) K_h(U_k - U_j)$$
$$\times \left(1 + \mu_2^{-1}\frac{U_i - U_j}{h}\frac{U_k - U_j}{h}\right)\frac{1}{f(U_j)}$$
$$\times \frac{\partial G(\varepsilon_i^*)^\tau}{\partial \varepsilon} V^{-1}(U_j) G(\varepsilon_k) X_i^\tau \Gamma^{-1}(U_j) X_k X_k^\tau A(U_k),$$

$$W_{n3} = \frac{1}{2n^2} \sum_{i=1}^n \sum_{k=1}^n \sum_{j=1}^n K_h(U_i - U_j) K_h(U_k - U_j)$$
$$\times \left(1 + \mu_2^{-1}\frac{U_i - U_j}{h}\frac{U_k - U_j}{h}\right)\frac{1}{f(U_j)}$$
$$\times \frac{\partial G(\varepsilon_i^*)^\tau}{\partial \varepsilon} V^{-1}(U_j)$$
$$\times \frac{\partial G(\varepsilon_k^*)}{\partial \varepsilon} A(U_i)^\tau X_i X_i^\tau \Gamma^{-1}(U_j) X_k X_k^\tau A(U_k),$$

where $\varepsilon_i^*$ is between $\varepsilon_i$ and $\varepsilon_i + A(U_i)^\tau X_i$ and $\varepsilon_k^*$ is between $\varepsilon_k$ and $\varepsilon_k + A(U_k)^\tau X_k$. Under some regularity conditions,

$$W_{n1} = \frac{1}{2n^2} \sum_{i=1}^n \sum_{k=1}^n G(\varepsilon_i)^\tau \Bigg\{ \sum_{j=1}^n K_h(U_i - U_j) K_h(U_k - U_j)$$
$$\times \left(1 + \mu_2^{-1} \frac{U_i - U_j}{h} \frac{U_k - U_j}{h}\right)$$
$$\times \frac{1}{f(U_j)} V^{-1}(U_j)$$
$$\times X_i^\tau \Gamma^{-1}(U_j) X_k X_k^\tau A(U_k) \Bigg\} \frac{\partial G(\varepsilon_k^*)}{\partial \varepsilon}$$
$$= \frac{W_{1n}^*}{2} + o_p(h^{-1/2}),$$
$$W_{n2} = W_{n1},$$

where $W_{1n}^*$ is defined in (3.5). Similarly, we write

$$W_{n3} = W_{n31} + 2W_{n32} + W_{n33}$$

where, when $EA(U)^\tau X X^\tau A(U) = O(\frac{1}{nh})$,

$$W_{n31} = \frac{1}{2n^2} \sum_{i=1}^n \sum_{k=1}^n \Xi_i^\tau \sum_{j=1}^n K_h(U_i - U_j) K_h(U_k - U_j)$$
$$\times \left(1 + \mu_2^{-1} \frac{U_i - U_j}{h} \frac{U_k - U_j}{h}\right) \frac{1}{f(U_j)} V^{-1}(U_j)$$
$$\times \Xi_k A(U_i)^\tau X_i X_i^\tau \Gamma^{-1}(U_j) X_k X_k^\tau A(U_k)$$
$$= \frac{1}{2n} \sum_{i=1}^n \sum_{k=1}^n \Xi_i^\tau K_h^*(U_i - U_k) V^{-1}(U_i)$$
$$\times \Xi_i A(U_i)^\tau X_i X_i^\tau \Gamma^{-1}(U_j) X_k X_k^\tau A(U_k) + o_p(h^{-1/2})$$
$$= O\left(\frac{1}{nh^2}\right) + \frac{W_{2n}^*}{2} + o_p(h^{-1/2}),$$
$$W_{n32} = \frac{1}{2n^2} \sum_{i=1}^n \sum_{k=1}^n \Xi_i^\tau \sum_{j=1}^n K_h(U_i - U_j) K_h(U_k - U_j)$$
$$\times \left(1 + \mu_2^{-1} \frac{U_i - U_j}{h} \frac{U_k - U_j}{h}\right) \frac{1}{f(U_j)} V^{-1}(U_j)$$
$$\times E\left[\frac{\partial G(\varepsilon_k^*)}{\partial \varepsilon} \bigg| U_k\right] A(U_i)^\tau X_i X_i^\tau \Gamma^{-1}(U_j) X_k X_k^\tau A(U_k)$$

$$= \frac{W_{3n}^*}{2} + o_p(h^{-1/2}),$$

$$W_{n33} = \frac{1}{2n^2} \sum_{i=1}^{n} \sum_{k=1}^{n} E\left[\frac{\partial G(\varepsilon_i^*)}{\partial \varepsilon}\Big|U_i\right] \sum_{j=1}^{n} K_h(U_i - U_j) K_h(U_k - U_j)$$

$$\times \left(1 + \mu_2^{-1} \frac{U_i - U_j}{h} \frac{U_k - U_j}{h}\right)$$

$$\times \frac{1}{f(U_j)} V^{-1}(U_j) E\left[\frac{\partial G(\varepsilon_k^*)}{\partial \varepsilon}\Big|U_k\right]$$

$$\times A(U_i)^\tau X_i X_i^\tau \Gamma^{-1}(U_j) X_k X_k^\tau A(U_k)$$

$$= O_p\left(\frac{1}{nh^2}\right)$$

$$+ \frac{n}{2} E\left\{E\left[\frac{\partial G(\varepsilon)}{\partial \varepsilon}\Big|U\right]^\tau V^{-1}(U) E\left[\frac{\partial G(\varepsilon)}{\partial \varepsilon}\Big|U\right] A(U)^\tau X X^\tau A(U)\right\}$$

$$\times (1 + o(1)),$$

with $\Xi_i$ defined in (3.6). Recall that $W_{2n}^*$ and $W_{3n}^*$ are in (3.7) and (3.8), respectively. Observe that as $EA(U)^\tau X X^\tau A(U) = O(\frac{1}{nh})$ we have

$$W_{n31} = \frac{1}{2n} \sum_{i=1}^{n} \sum_{k=1}^{n} \Xi_i^\tau K_h^*(U_i - U_k) V^{-1}(U_i)$$

$$\times \Xi_i A(U_i)^\tau X_i X_i^\tau \Gamma^{-1}(U_j) X_k X_k^\tau A(U_k) + o_p(h^{-1/2})$$

$$= O\left(\frac{1}{nh^2}\right) + \frac{W_{2n}^*}{2} + o_p(h^{-1/2}),$$

$$W_{n32} = W_{3n}^*/2 + o_p(h^{-1/2}),$$

$$W_{n33} = O_p\left(\frac{1}{nh^2}\right)$$

$$+ \frac{n}{2} E\left\{E\left[\frac{\partial G(\varepsilon)}{\partial \varepsilon}\Big|U\right]^\tau V^{-1}(U) E\left[\frac{\partial G(\varepsilon)}{\partial \varepsilon}\Big|U\right] A(U)^\tau X X^\tau A(U)\right\}$$

$$\times (1 + o(1)).$$

Similarly we have

$$l(G) = (1 + o_p(h^{-1/2}))$$

$$\times \left[\frac{1}{2n^2} \sum_{i=1}^{n} \sum_{k=1}^{n} \sum_{j=1}^{n} K_h(U_i - U_j) K_h(U_k - U_j)\right.$$

(5.13)

$$\times \left(1 + \mu_2^{-1} \frac{U_i - U_j}{h} \frac{U_k - U_j}{h}\right)$$

$$\times \frac{1}{f(U_j)} G^\tau(\varepsilon_i) C(U_j) G(\varepsilon_k) + 2 S_{n1} + S_{n2}\bigg]$$

where

$$S_{n1} = \frac{1}{2n^2} \sum_{i=1}^{n} \sum_{k=1}^{n} \sum_{j=1}^{n} K_h(U_i - U_j) K_h(U_k - U_j)$$

$$\times \left(1 + \mu_2^{-1} \frac{U_i - U_j}{h} \frac{U_k - U_j}{h}\right) \frac{1}{f(U_j)}$$

$$\times G^\tau(\varepsilon_i) C(U_j) \frac{\partial G(\varepsilon_k^*)}{\partial \varepsilon} A''(U_j^*)^\tau (U_k - U_j)^2 X_k,$$

$$= O_p(n(nh)^{-1} h^2)$$

$$= O_p(h)$$

and

$$S_{n2} = \frac{1}{2n^2} \sum_{i=1}^{n} \sum_{k=1}^{n} \sum_{j=1}^{n} K_h(U_i - U_j) K_h(U_k - U_j)$$

$$\times \left(1 + \mu_2^{-1} \frac{U_i - U_j}{h} \frac{U_k - U_j}{h}\right) \frac{1}{f(U_j)}$$

$$\times \frac{\partial G(\varepsilon_i^*)^\tau}{\partial \varepsilon} C(U_j) \frac{\partial G(\varepsilon_k^*)}{\partial \varepsilon} A''(U_j) X_i X_i^\tau \Gamma^{-1}(U_j)$$

$$\times X_k X_k^\tau A''(U_j)(U_i - U_j)^2 (U_k - U_j)^2$$

$$= \frac{nh^4}{8} E\{D^\tau(U) C(U) D(U) A''(U)^\tau X X^\tau A''(U)\}$$

$$\times \iint t^2 (s+t)^2 K(t) K(s+t)(1 + \mu_2^{-1} t(s+t)) \, dt \, ds \, (1 + o_p(1)),$$

where $U_j^*$ is between $U_k$ and $U_j$. Now the desired result follows from (5.12) and (5.13). This proves the theorem. $\square$

## APPENDIX

LEMMA 1. *Under conditions* (K0), (U0), (A2)–(A4), *if there exist some positive constants* $b_0, b_1$ *and* $\eta < 1/2$ *such that* $b_0 \leq hn^\eta \leq b_1$, *then there exists a sequence of positive constants* $d_n \to 0$ *such that*

$$A_n(u_0, \beta) = E\left\{K_h(U - u_0) G\left(Y - \beta^\tau Z\left(X, \frac{U - u_0}{h}\right)\right) \otimes Z\left(X, \frac{U - u_0}{h}\right)\right\}$$

$$+ o_p(n^{-1/\xi} \wedge h^{1/2}) d_n.$$

*Furthermore, if condition* (A5) *holds and* $\eta > 1/(2\xi)$, *then uniformly in* $\|\beta - \beta_0\| \leq r_n = o(n^{-1/\xi})d_n$,

$$A_n(u_0, \beta) = o_p(n^{-1/\xi})d_n.$$

PROOF.  For any positive constant $M_n$, we can write

$$A_n(u_0, \beta) = EK_h(U - u_0)G\left(Y - \beta^\tau Z\left(X, \frac{U - u_0}{h}\right)\right) \otimes Z\left(X, \frac{U - u_0}{h}\right)$$
$$+ A_{n1}(u_0, \beta) + A_{n2}(u_0, \beta),$$

where

$$A_{n1}(u_0, \beta) = \frac{1}{n}\sum_{i=1}^n K_h(U_i - u_0)\mathbf{G}_{ih}(u_0, \beta)I(F(Y_i, X_i) \leq M_n)$$
$$- EK_h(U - u_0)\mathbf{G}_h(u_0, \beta)I(F(Y, X) \leq M_n)$$

and

$$A_{n2}(u_0, \beta) = \frac{1}{n}\sum_{i=1}^n K_h(U_i - u_0)\mathbf{G}_{ih}(u_0, \beta)I(F(Y_i, X_i) > M_n)$$
$$- EK_h(U - u_0)\mathbf{G}_h(u_0, \beta)I(F(Y, X) > M_n).$$

Note that

(A.1) $$E\|A_{n2}(u_0, \beta)\| \leq 2EK_h(U - u_0)\mathbf{G}_h(u_0, \beta)I(F(Y, X) > M_n) \leq cM_n^{1-\xi}.$$

Consider the following empirical processes:

$$v_n(g) = n^{-1/2}\sum_{i=1}^n (g(Y_i, X_i, u_0, \beta) - Eg(Y, X, u_0, \beta)),$$
$$g \in \mathcal{F}_n = \{M_n^{-1}g : g \in \mathcal{F}_1\},$$

where $\mathcal{F}_1$ is defined as in Section 5.1. It follows directly from assumption (A4) that

$$N(\delta, L_1(P_n), \mathcal{F}_n) \leq c_1(P_n)(h\delta M_n)^{-w_1}.$$

Obviously, by condition (A3), for $g = K((\cdot - u_0)/h)\mathbf{G}(u_0, \beta) \in \mathcal{F}_n$,

$$E\|g(Y, X, u_0, \beta)\|^2$$
$$\leq chM_n^{-2} \sup_{u_0, t, \|\beta - \beta_0\| \leq r_0} E_{Y|U = u_0 + th}\{G_k^2(Y - \beta^\tau Z(X, t))\|Z(X, t)\|^2\}$$
$$\leq O(hM_n^{-2}) = v.$$

Now let $M_n = n^{s_0}$, $\delta_n = (h^{1/2} \wedge n^{-1/\xi})(\log n)^{-1}$ and $M = \delta_n n^{1/2} h M_0 M_n^{-1}$. Using Lemma 2 in Zhang and Gijbels (2003), we have

(A.2)
$$\begin{aligned}
&P\{\sup \|A_{n1}(u_0,\beta)\|\delta_n^{-1} > M_0\} \\
&= P\Big\{\sup_{g \in \mathcal{F}_n} \|v_n(g)\| > M\Big\} \\
&\leq c_1(n^{1/2}(MhM_n)^{-1})^{w_1} \exp\{-c_3 M^2/v\} + c_2 v^{-w_1} \exp(-nv) \\
&= O((h^2\delta_n)^{-w_1}) \exp\{-c_3 \delta_n^2 n h^2 M_0^2 M_n^{-2}/hM_n^{-2}\} \\
&\quad + c_2 O(hM_n^{-2})^{-w_1} \exp(-c_4 n h M_n^{-2}).
\end{aligned}$$

The last terms in (A.1) and (A.2) are $o(\delta_n)$ and $o(1)$, respectively, if

$$b_0 \leq h n^\eta \leq b_1, \qquad nh^2/\log n \to \infty, \qquad n^{1-2/\xi}h/\log n \to \infty,$$
$$nhM_n^{-2}/\log n \to \infty, \qquad M_n^{-\xi+1}\delta_n^{-1} \to 0.$$

The above requirements are fulfilled provided that, for $s_0 > 0$,

$$b_0 \leq h n^\eta \leq b_1, \qquad 0 < \eta < \min\Big\{\frac{1}{2}, 1 - \frac{2}{\xi}\Big\},$$
$$\max\Big\{\frac{\eta}{2(\xi-1)}, \frac{1}{\xi(\xi-1)}\Big\} < s_0 < \frac{1-\eta}{2}.$$

These conditions are equivalent to

$$b_0 \leq h n^\eta \leq b_1, \qquad 0 < \eta < \min\Big\{\frac{1}{2}, 1 - \frac{2}{\xi}, 1 - \frac{1}{\xi}, 1 - \frac{2}{\xi(\xi-1)}\Big\} = \frac{1}{2}$$

since $\xi \geq 4$.

Let $d_n = (\log n)^{-1}$ and $f$ be the density of $U$. Now we can complete the proof if we note that for $\|\beta - \beta_0\| \leq o(n^{-1/\xi})d_n$ and $\eta > 1/(2\xi)$, we have

$$\begin{aligned}
&E\Big\{K_h(U - u_0)EK_h(U - u_0)G\Big(Y - \beta^\tau Z\Big(X, \frac{U-u_0}{h}\Big)\Big) \otimes Z\Big(X, \frac{U-u_0}{h}\Big)\Big\} \\
&= \int K(t)\Big\{E\Big[G\Big(Y - \beta^\tau Z\Big(X, \frac{U-u_0}{h}\Big)\Big)\Big|U = u_0 + th\Big]\Big\} \\
&\quad \otimes E[Z(X,t)|U = u_0 + th]f(u_0+th)\,dt \\
&= O(h^2) + O(\|\beta - \beta_0\|) \\
&= O(n^{-2\eta}) + o(n^{-1/\xi}d_n) \\
&= o(n^{-1/\xi}d_n)
\end{aligned}$$

by using condition (A5). □

LEMMA 2. *Under conditions* (K0), (U0), (A2), (A6) *and* (A7), *as* $n \to \infty$, $b_0 \leq hn^\eta \leq b_1$, $0 < \eta < 1/2$, *we have*

$$V_n(u_0, \beta) = EK_h(U - u_0)\mathbf{G}_h(u_0, \beta)\mathbf{G}_h^\tau(u_0, \beta) + o_p(h^{1/2})$$
$$= V(u_0) \otimes (S \otimes \Gamma(u_0)) + o_p(h^{1/2}) + O(\|\beta - \beta_0\|).$$

PROOF. The proof is similar to that of Lemma 1 and is thus omitted. □

LEMMA 3. *Under conditions* (K0), (U0), (A1)–(A10) *and* (B1), *if* $b_0 \leq hn^\eta \leq b_1$, $1/(2\xi) < \eta < 1/2$, *then there exists a sequence of positive constants* $d_n \to 0$ *such that as* $n \to \infty$,

$$\hat{\beta}(u_0) = \beta_0(u_0) + o_p(n^{-1/\xi} \wedge h^{1/2})d_n,$$
$$\alpha_n(u_0, \hat{\beta}(u_0)) = o_p(n^{-1/\xi} \wedge h^{1/2}).$$

PROOF. First of all, by Lemma 1, there exists a sequence of positive constants $d_n \to 0$ such that

(A.3) $$A_n(u_0, \beta_0) = o_p(n^{-1/\xi} \wedge h^{1/2})d_n.$$

Note that condition (A2) implies

(A.4) $$Z_n(u_0, \beta) = o_p(n^{1/\xi})$$

uniformly in $u_0 \in \Omega$ and $\|\beta - \beta_0\| \leq r_0$. Set the function

$$g_n(\alpha, \beta) = \frac{1}{n}\sum_{i=1}^n K_h(U_i - u_0)\frac{\mathbf{G}_{ih}(u_0, \beta)}{1 + \alpha^\tau \mathbf{G}_{ih}(u_0, \beta)}.$$

Then following the argument of Owen (1990) and using conditions (K0), (U0), (A1), (A4), (A5), (A8), (A9) and (B1), we can show that for large $n$, $\alpha_n(u_0, \beta)$ exists and satisfies the equation

(A.5) $$g_n(\alpha_n(u_0, \beta), \beta) = 0$$

when $\|\beta - \beta_0\| \leq r_0$ and $r_0$ is small. To see this, we first note that for constant $\delta > 0$ small enough, we have

$$\inf_{\substack{\|\psi\|=1 \\ u_0 \in \Omega}} \int K(t) E[I\{\psi^\tau(G(\varepsilon) \otimes Z(X,t)) > \delta\}|U = u_0]\, dt > \delta,$$

which yields

(A.6) $$\inf_{\substack{\|\psi\|=1 \\ u_0 \in \Omega \\ \|\beta-\beta_0\| \leq r_0}} \int K(t) E[I\{\psi^\tau \mathbf{G}_h(u_0, \beta) > \delta/2\}|U = u_0 + th]$$
$$\times f(u_0 + th)\, dt \geq \delta/2$$

as $h \to 0$ and $r_0$ is small enough. This is the main consequence of conditions (A1) and (A9). Define

$$H_n(\beta, \psi) = \frac{1}{n} \sum_{i=1}^n w_h(U_i, u_0) I\{\mathbf{G}_{ih}(u_0, \beta)^\tau \psi > \delta\}.$$

Then under conditions (A1)–(A4), (A8)–(A10), using (A.6) and the strong convergence of empirical processes [Pollard (1984) and van der Vaart and Wellner (1996)], we can show that there exists $\delta > 0$ such that for small $r_0$ and large $n$, $\inf_{\|\psi\|=1} H_n(\beta, \psi) > \delta$ almost surely. This shows that 0 is contained in the convex hull of the points in $\{\mathbf{G}_{ih}(u_0, \beta) : w_h(U_i, u_0) > 0, 1 \leq i \leq n\}$. Now (A.5) follows directly from the Lagrange multiplier method as in Owen (1990).

Let

$$\alpha_n(u_0, \beta_0) = \rho v \qquad \text{with } \rho = \|\alpha_n(u_0, \beta_0)\| \text{ and } \|v\| = 1.$$

We have

$$\begin{aligned}
0 &= \|g_n(\alpha_n(u_0, \beta_0), \beta_0)\| \\
&= \|g_n(\rho v, \beta_0)\| \\
&\geq \|v^\tau g_n(\rho v, \beta_0)\| \\
&= \frac{1}{n} \left| v^\tau \sum_{i=1}^n K_h(U_i - u_0) \left\{ \mathbf{G}_{ih}(u_0, \beta_0) - \frac{\rho \mathbf{G}_{ih}(u_0, \beta_0) \mathbf{G}_{ih}(u_0, \beta_0)^\tau v}{1 + \rho v^{*\tau} \mathbf{G}_{ih}(u_0, \beta_0)} \right\} \right| \\
&\geq \frac{1}{n} \rho \sum_{i=1}^n K_h(U_i - u_0) \frac{v^\tau \mathbf{G}_{ih}(u_0, \beta_0) \mathbf{G}_{ih}(u_0, \beta_0)^\tau v}{1 + \rho v^{*\tau} \mathbf{G}_{ih}(u_0, \beta_0)} - |v^\tau A_n(u_0, \beta_0)| \\
&\geq \rho \frac{v^\tau V_n(u_0, \beta_0) v}{1 + \rho Z_n(u_0, \beta_0)} - \|A_n(u_0, \beta_0)\|,
\end{aligned}$$

where $v^* = tv$ with $0 \leq t \leq 1$. Thus, combining (A.4) with (A.3), Lemma 2 and condition (A9), we have

$$\begin{aligned}
\rho &\leq \frac{\|A_n(u_0, \beta_0)\|}{v^\tau V_n(u_0, \beta_0) v - \|A_n(u_0, \beta_0)\| Z_n(u_0, \beta_0)} \\
&= O_p(\|A_n(u_0, \beta_0)\|) \\
&= o_p(n^{-1/\xi} \wedge h^{1/2}) d_n,
\end{aligned}$$

that is,

(A.7) $$\alpha_n(u_0, \beta_0) = o_p(n^{-1/\xi} \wedge h^{1/2}) d_n.$$

Set $\phi_n = (h^{1/2} \wedge n^{-1/\xi}) d_n$, and let $u(u_0, \beta)$ satisfy

$$u(u_0, \beta) \| E\{K_h(U - u_0) \mathbf{G}_h(u_0, \beta)\} \| = E\{K_h(U - u_0) \mathbf{G}_h(u_0, \beta)\}.$$

Define

$$l_n(u_0, \beta) = -\frac{1}{n} \sum_{i=1}^n K_h(U_i - u_0) \log(1 + \alpha_n(u_0, \beta)^\tau \mathbf{G}_{ih}(u_0, \beta)),$$

$$T_n(u_0, \beta) = \frac{1}{n} \sum_{i=1}^n K_h(U_i - u_0) \log(1 + \phi_n u(u_0, \beta)^\tau \mathbf{G}_{ih}(u_0, \beta)),$$

$$T_{n1}(u_0, \beta) = \frac{1}{n} \sum_{i=1}^n K_h(U_i - u_0) \log(1 + \phi_n u(u_0, \beta)^\tau \mathbf{G}_{ih}(u_0, \beta))$$
$$\times I(\|\mathbf{G}_{ih}\| \leq n^{1/\xi}).$$

We have

(A.8) $\qquad 0 \geq l_n(u_0, \beta_0) \geq -\alpha_n^\tau(u_0, \beta) A_n(u_0, \beta) = o_p(\phi_n^2),$

and uniformly for $u_0$ and $\beta$,

$$T_{n1}(u_0, \beta) = \phi_n \frac{1}{n} \sum_{i=1}^n K_h(U_i - u_0) u(u_0, \beta)^\tau \mathbf{G}_{ih}(u_0, \beta) I(\|\mathbf{G}_{ih}\| \leq n^{1/\xi})$$
$$- \frac{1}{2} \phi_n^2 |O(1)| \frac{1}{n} \sum_{i=1}^n K_h(U_i - u_0) F(Y_i, X_i)^2$$
$$= \phi_n \{u(u_0, \beta)^\tau E[K_h(U - u_0) \mathbf{G}_h(u_0, \beta)]\} + o_p(\phi_n^2) + O_p(\phi_n^2).$$

Note that for fixed $u_0$ and $\beta$, the function $-\frac{1}{n} \sum_{i=1}^n K_h(U_i - u_0) \log(1 + \alpha^\tau \mathbf{G}_{ih}(u_0, \beta))$ attains the minimum at $\alpha_n(u_0, \beta)$. This implies $l_n(u_0, \beta) \leq -T_n(u_0, \beta)$. Consequently, for any $\rho > 0$, by (A.8) we have

$$P(\|\hat{\beta}(u_0) - \beta_0\| > \rho)$$
$$\leq P\left(\sup_{\|\beta - \beta_0\| \geq \rho} l_n(u_0, \beta) > l_n(u_0, \beta_0) \text{ for some } u_0\right)$$
$$\leq P\left(\sup_{\|\beta - \beta_0\| \geq \rho} l_n(u_0, \beta) > -|O_p(\phi_n^2)| \text{ for some } u_0\right)$$
$$\leq P\left(\sup_{\|\beta - \beta_0\| \geq \rho} (-T_n(u_0, \beta)) \geq -|O_p(\phi_n^2)| \text{ for some } u_0\right)$$
$$\leq P\left(\sup_{\|\beta - \beta_0\| \geq \rho} (-T_{n1}(u_0, \beta)) \geq -|O_p(\phi_n^2)| \text{ for some } u_0\right)$$
$$+ P\left(\sup_{u_0, \beta} Z_n(u_0, \beta) > n^{1/\xi}\right)$$
$$\leq P\left\{\inf_{\|\beta - \beta_0\| \geq \rho} \|E K_h(U - u_0) G_h(u_0, \beta)\| \leq |O_p(\phi_n)| \text{ for some } u_0\right\} + o(1)$$

$$\to 0,$$

where the last limit follows from condition (A10). Therefore using for $\rho_{n1} \to 0$ and $\rho_2$ in condition (A10), as $n \to \infty$, we have

$$P(\|\hat{\beta}(u_0) - \beta_0\| > \rho_{n1} \text{ for some } u_0)$$
$$= P(\rho_2 \geq \|\hat{\beta}(u_0) - \beta_0\| > \rho_{n1} \text{ for some } u_0) + o(1)$$
$$\leq P\bigg(\inf_{\rho_2 \geq \|\beta - \beta_0\| \geq \rho_{n1}} \|EK_h(U - u_0)G_h(u_0, \beta)\| \leq |O_p(\phi_n)| \text{ for some } u_0\bigg)$$
$$+ o(1)$$
$$\leq P(\rho_{n1} + O(h^2) \leq |O_p(\phi_n)|) + o(1),$$

which leads to

$$\hat{\beta}(u_0) - \beta_0 = O_p(\phi_n) = o_p(n^{-1/\xi} \wedge h^{1/2})d_n.$$

Invoking the argument of Owen (1990) and Lemma 1 again, we have

$$\alpha_n(u_0, \hat{\beta}(u_0)) = o_p(n^{-1/\xi} \wedge h^{1/2})$$

uniformly in $u_0$. This completes the proof. □

LEMMA 4. *Suppose for some positive constants $b_0$ and $b_1$, $b_0 \leq hn^\eta \leq b_1$, $0 < \eta < 1/2$. Then under conditions* (K0), (U0), (A2), (A6), (A7) *and* (A9), *as $n \to \infty$, we have*

$$V_n(u_0, \alpha, \beta) = V(u_0) \otimes (S \otimes \Gamma(u_0))(1 + o_p(h^{1/2}))$$

*uniformly for $u_0 \in \Omega$, $\|\alpha\| + \|\beta - \beta_0\| \leq o(n^{-1/\xi} \wedge h^{1/2})$.*

PROOF. Note that under condition (A2) we have

$$\sup_{\substack{u_0 \in \Omega \\ \|\beta - \beta_0\| \leq r_0}} Z_n(u_0, \beta) = o_p(n^{1/\xi}),$$

which together with Lemma 2 yields

$$V_n(u_0, \alpha, \beta) = V_n(u_0, \beta) + O_p(\|\alpha\|)\frac{(2 + o_p(1))}{(1 + o_p(1))}\frac{1}{n}\sum_{i=1}^n K_h(U_i - u_0)F(Y_i, X_i)^3$$
$$= V_n(u_0, \psi_2) + O_p(\|\alpha\|)$$
$$= V(u_0) \otimes (S \otimes \Gamma(u_0)) + o_p(h^{1/2}) + O_p(\|\alpha\|).$$

The proof is complete. □

LEMMA 5. *Suppose there exist positive constants $b_0$, $b_1$ and $\eta$ such that $b_0 \leq hn^\eta \leq b_1$, $0 < \eta < 1/2$. Then under conditions* (K0), (U0), (A2), (B1)–(B4), *as $n \to \infty$,*
$$B_n(u_0, \alpha, \beta) = D(u_0) \otimes (S \otimes \Gamma(u_0))(1 + o_p(h^{1/2}))$$
*uniformly for $u_0 \in \Omega$, $\|\alpha\| + \|\beta - \beta_0\| \leq o(n^{-1/\xi} \wedge h^{1/2})$.*

The proof is similar to that of Lemma 1 and thus is omitted.

LEMMA 6. *Under conditions* (K0), (U0), (A2), (B1), *as $h \to 0$, $nh \to \infty$,*
$$C_n(u_0, \alpha, \beta) = O_p(\|\alpha\|)$$
*uniformly for $u_0 \in \Omega$, $\|\alpha\| + \|\beta - \beta_0\| \leq o(n^{-1/\xi} \wedge h^{1/2})$.*

PROOF. Note that by condition (A2) and $\|\psi_1\| \leq o(n^{-1/\xi} \wedge h^{1/2})$, we have
$$\max_i \sup_{\beta, u_0} \|\alpha^\tau \mathbf{G}_{ih}(u_0, \beta)\| = o_p(1).$$
Thus
$$\|C_n(u_0, \alpha, \beta)\| \leq O_p(\|\alpha\|) \frac{1}{n} \sum_{i=1}^n K_h(U_i - u_0) F_4(Y_i, X_i) F(Y_i, X_i) = O_p(\|\alpha\|)$$
by conditions (A2) and (B1). The proof is complete. □

LEMMA 7. *Under conditions* (K0), (U0) *and* (B5), *as $h \to 0$ and $nh \to \infty$,*
$$D_n(u_0, \alpha, \beta) = O_p(\|\alpha\|)$$
*uniformly for $u_0 \in \Omega$, $\|\alpha\| + \|\beta - \beta_0\| \leq o(n^{-1/\xi} \wedge h^{1/2})$.*

LEMMA 8. *Under conditions* (K0), (U0) *and* (B5), *as $h \to 0$, $nh \to \infty$,*
$$E_n(u_0, \alpha, \beta) = O_p(\|\alpha\|^2)$$
*uniformly for $u_0 \in \Omega$, $\|\alpha\| + \|\beta - \beta_0\| \leq o(n^{-1/\xi} \wedge h^{1/2})$.*

The proofs of Lemmas 7 and 8 are similar to the proof of Lemma 6 and thus are omitted.

**Acknowledgments.** The authors thank the Editor, Jon A. Wellner, an Associate Editor and two anonymous referees for their constructive comments that have helped us improve the presentation and results of the paper. The work was partially done when the second author was visiting Wing Hung Wong's laboratory. Wing Hung Wong's support is greatly appreciated.

# REFERENCES

Azzalini, A. and Bowman, A. N. (1993). On the use of nonparametric regression for checking linear relationships. *J. Roy. Statist. Soc. Ser. B* **55** 549–557. MR1224417

Bickel, P. J., Klaassen, C. A. J., Ritov, Y. and Wellner, J. (1993). *Efficient and Adaptive Estimation for Semiparametric Models.* Johns Hopkins Univ. Press, Baltimore, MD. MR1245941

Bickel, P. J. and Ritov, Y. (1992). Testing for goodness of fit: A new approach. In *Nonparametric Statistics and Related Topics* (A. K. Md. E. Saleh, ed.) 51–57. North-Holland, Amsterdam. MR1226715

Brumback, B. and Rice, J. A. (1998). Smoothing spline models for the analysis of nested and crossed samples of curves (with discussion). *J. Amer. Statist. Assoc.* **93** 961–994. MR1649194

Cai, Z., Fan, J. and Li, R. (2000). Efficient estimation and inferences for varying-coefficient models. *J. Amer. Statist. Assoc.* **95** 888–902. MR1804446

Cai, Z., Fan, J. and Yao, Q. (2000). Functional-coefficient regression models for nonlinear time series. *J. Amer. Statist. Assoc.* **95** 941–956. MR1804449

Carroll, R. J., Ruppert, D. and Welsh, A. H. (1998). Local estimating equations. *J. Amer. Statist. Assoc.* **93** 214–227. MR1614624

Cleveland, W. S., Grosse, E. and Shyu, W. M. (1991). Local regression models. In *Statistical Models in S* (J. M. Chambers and T. J. Hastie, eds.) 309–376. CRC Press, Boca Raton, FL.

Chen, R. and Tsay, R. S. (1993). Functional-coefficient autoregressive models. *J. Amer. Statist. Assoc.* **88** 298–308. MR1212492

Chen, S. X., Gao, J. and Li, M. (2003). Simultaneous specification tests for the mean and variance structures of regression with applications to testing of diffusion models. Unpublished manuscript.

Chen, S. X., Härdle, W. and Li, M. (2003). An empirical likelihood goodness-of-fit test for time series. *J. R. Stat. Soc. Ser. B Stat. Methodol.* **65** 663–678. MR1998627

Eubank, R. L. and Hart, J. D. (1992). Testing goodness-of-fit in regression via order selection criteria. *Ann. Statist.* **20** 1412–1425. MR1186256

Eubank, R. L. and LaRiccia, V. N. (1992). Asymptotic comparison of Cramér–von Mises and nonparametric function estimation techniques for testing goodness-of-fit. *Ann. Statist.* **20** 2071–2086. MR1193326

Fan, J. (1992). Design-adaptive nonparametric regression. *J. Amer. Statist. Assoc.* **87** 998–1004. MR1209561

Fan, J. (1996). Test of significance based on wavelet thresholding and Neyman's truncation. *J. Amer. Statist. Assoc.* **91** 674–688. MR1395735

Fan, J. and Huang, L. (2001). Goodness-of-fit test for parametric regression models. *J. Amer. Statist. Assoc.* **96** 640–652. MR1946431

Fan, J., Zhang, C. and Zhang, J. (2001). Generalized likelihood ratio statistics and Wilks phenomenon. *Ann. Statist.* **29** 153–193. MR1833962

Fan, J. and Zhang, W. (1999). Statistical estimation in varying coefficient models. *Ann. Statist.* **27** 1491–1518. MR1742497

Fan, Y. and Li, Q. (1996). Consistent model specification tests: Omitted variables and semiparametric functional forms. *Econometrica* **64** 865–890. MR1399221

Farrell, R. H. (1985). *Multivariate Calculation. Use of the Continuous Groups.* Springer, New York. MR770934

Härdle, W. and Mammen, E. (1993). Comparing nonparametric versus parametric regression fits. *Ann. Statist.* **21** 1926–1947. MR1245774

Hart, J. D. (1997). *Nonparametric Smoothing and Lack-of-Fit Tests*. Springer, New York. MR1461272

Hastie, T. J. and Tibshirani, R. J. (1993). Varying-coefficient models (with discussion). *J. Roy. Statist. Soc. Ser. B* **55** 757–796. MR1229881

Hettmansperger, T. P. (1984). *Statistical Inference Based on Ranks.* Wiley, New York. MR758442

Hong, Y. and Lee, T.-H. (2003). Inference on predictability of foreign exchange rates via generalized spectrum and nonlinear time series models. *Rev. Econom. Statist.* **85** 1048–1062.

Horowitz, J. L. and Spokoiny, G. G. (2001). An adaptive, rate-optimal test of a parametric mean-regression model against a nonparametric alternative. *Econometrica* **69** 599–631. MR1828537

Horowitz, J. L. and Spokoiny, G. G. (2002). An adaptive, rate-optimal test of linearity for median regression models. *J. Amer. Statist. Assoc.* **97** 822–835. MR1941412

Inglot, T. and Ledwina, T. (1996). Asymptotic optimality of data-driven Neyman's tests for uniformity. *Ann. Statist.* **24** 1982–2019. MR1421157

Ingster, Yu. I. (1993). Asymptotically minimax hypothesis testing for nonparametric alternatives. I–III. *Math. Methods Statist.* **2** 85–114, 171–189, 249–268. MR1257978

Kallenberg, W. C. M. and Ledwina, T. (1997). Data-driven smooth tests when the hypothesis is composite. *J. Amer. Statist. Assoc.* **92** 1094–1104. MR1482140

Kitamura, Y. (1997). Empirical likelihood methods with weakly dependent processes. *Ann. Statist.* **25** 2084–2102. MR1474084

LeBlanc, M. and Crowley, J. (1995). Semiparametric regression functionals. *J. Amer. Statist. Assoc.* **90** 95–105. MR1325117

Newey, W. K. (1993). Efficient estimation of models with conditional moment restrictions. In *Econometrics. Handbook of Statistics* **11** (G. S. Maddala, C. R. Rao and H. D. Vinod, eds.) 419–453. North-Holland, Amsterdam. MR1247253

Owen, A. B. (1988). Empirical likelihood ratio confidence intervals for a single functional. *Biometrika* **75** 237–249. MR946049

Owen, A. B. (1990). Empirical likelihood ratio confidence regions. *Ann. Statist.* **18** 90–120. MR1041387

Pollard, D. (1984). *Convergence of Stochastic Processes*. Springer, New York. MR762984

Qin, J. and Lawless, J. (1994). Empirical likelihood and general estimating equations. *Ann. Statist.* **22** 300–325. MR1272085

Shiryayev, A. N. (1996). *Probability Theory*, 2nd ed. Springer, New York.

Spokoiny, V. G. (1996). Adaptive hypothesis testing using wavelets. *Ann. Statist.* **24** 2477–2498. MR1425962

Tong, H. (1990). *Nonlinear Time Series*: *A Dynamical System Approach*. Oxford Univ. Press. MR1079320

van der Vaart, A. W. and Wellner, J. A. (1996). *Weak Convergence and Empirical Processes.* Springer, New York. MR1385671

Wilks, S. S. (1938). The large-sample distribution of the likelihood ratio for testing composite hypotheses. *Ann. Math. Statist.* **9** 60–62.

Wu, C. O., Chiang, C.-T. and Hoover, D. R. (1998). Asymptotic confidence regions for kernel smoothing of a varying-coefficient model with longitudinal data. *J. Amer. Statist. Assoc.* **93** 1388–1401. MR1666635

Zhang, C. M. (2003). Adaptive tests of regression functions via multi-scale generalized likelihood ratios. *Canad. J. Statist.* **31** 151–171. MR2016225

Zhang, J. and Gijbels, I. (2003). Sieve empirical likelihood and extensions of generalized least squares. *Scand. J. Statist.* **30** 1–24. MR1963890

Zhang, J. and Liu, A. (2003). Local polynomial fitting based on empirical likelihood. *Bernoulli* **9** 579–605. MR1996271

Department of Operations Research
  and Financial Engineering
Princeton University
Princeton, New Jersey 08544
USA
e-mail: jqfan@princeton.edu

Institute of Mathematics and Statistics
University of Kent at Canterbury
Kent CT2 7NF
United Kingdom
e-mail: j.zhang@kent.ac.uk